\documentclass[seceqn,secthm,final]{elsart}
\usepackage{amssymb,amsmath,bm}

\begin{document}

\newcommand{\prs}{\textrm{\upshape prs}}
\newcommand{\rprs}{\textrm{\upshape rprs}}
\newcommand{\lc}{\textrm{\upshape lc}}
\newcommand{\syl}{\textrm{\upshape Syl}}
\newcommand{\subres}{\textrm{\upshape S}}
\newcommand{\recsubres}{\bar{\textrm{\upshape S}}} 
\newcommand{\nessubres}{\tilde{\textrm{\upshape S}}}
\newcommand{\rednessubres}{\hat{\textrm{\upshape S}}}

\journal{Journal of Algebra}
\begin{frontmatter}

\title{Recursive Polynomial Remainder Sequence and its
  Subresultants\thanksref{label1}}

\thanks[label1]{Preliminary versions of this paper have been presented
  at CASC 2003 (Passau, Germany, September 20--26, 2003)
  \cite{ter2003} and CASC 2005 (Kalamata, Greece, September
  12--16, 2005) \cite{ter2005}.}

\author{Akira Terui}
\address{Graduate School of Pure and Applied Sciences\\
  University of Tsukuba\\
  Tsukuba, 305-8571 Japan}
\ead{terui@math.tsukuba.ac.jp}
\date{January 25, 2008}

\begin{abstract}
  We introduce concepts of ``recursive polynomial remainder sequence
  (PRS)'' and ``recursive subresultant,'' along with investigation of
  their properties.  A \emph{recursive PRS} is defined as, if there
  exists the GCD (greatest common divisor) of initial polynomials, a
  sequence of PRSs calculated ``recursively'' for the GCD and its
  derivative until a constant is derived, and \emph{recursive
    subresultants} are defined by determinants representing the
  coefficients in recursive PRS as functions of coefficients of
  initial polynomials.  We give three different constructions of
  subresultant matrices for recursive subresultants; while the first
  one is built-up just with previously defined matrices thus the size
  of the matrix increases fast as the recursion deepens, the last one
  reduces the size of the matrix drastically by the Gaussian
  elimination on the second one which has a ``nested'' expression,
  \textit{i.e.\/} a Sylvester matrix whose elements are themselves
  determinants.
\end{abstract}

\end{frontmatter}

\section{Introduction}
\label{intro}

The polynomial remainder sequence (PRS) is one of the most fundamental
tools in computer algebra.  Although the Euclidean algorithm (see
Knuth \cite{knu1998}) for calculating PRS is simple, coefficient
growth in PRS makes the Euclidean algorithm often very inefficient.
To overcome this problem, the mechanism of coefficient growth has been
extensively studied through the theory of subresultants; see Collins
\cite{col1967}, Brown and Traub \cite{bro-tra71}, Loos
\cite{loos1983}, etc.  By the theory of subresultant, we can remove
extraneous factors of the elements of PRS systematically.

In this paper, we consider a variation of the subresultant.  When we
calculate PRS for polynomials which have a nontrivial GCD, we usually
stop the calculation with the GCD.  However, it is sometimes useful to
continue the calculation by calculating the PRS for the GCD and its
derivative; this is necessary for calculating the number of real zeros
including their multiplicities.  We call such a PRS a ``recursive
PRS.''

Although the theory of subresultants has been developed widely, the
corresponding theory for recursive PRS is still unknown within the
author's knowledge; this is the main problem which we investigate in
this paper.  By ``recursive subresultants,'' we denote determinants
which represent elements of recursive PRS as functions of the
coefficients of initial polynomials.

We give three different constructions of subresultant matrices to
express recursive subresultants in this paper.  The first matrix
construction recursively builds the matrix by shifting previously
defined matrices, similarly as the Sylvester matrix shifts
coefficients of the initial polynomials, thus the size of the matrices
increases fast as the recursion deepens.  The second matrix
construction uses ``nested'' matrices, or a Sylvester matrix whose
entries are themselves determinants.  Finally, by the Gaussian
elimination with the Sylvester's identity on the second construction,
we succeed to give the reduced matrix construction which expresses the
coefficients of the polynomials in the recursive PRS as determinants
of very small matrices, whose size actually decreases as the recursion
deepens.

This paper is organized as follows.  In Sect.~\ref{sec:recprs}, we
introduce the concept of recursive PRS.  In
Sect.~\ref{sec:subresrecprs}, we define recursive subresultant and
show its relationship to recursive PRS.  In Sect.~\ref{sec:nessubres},
we define the ``nested subresultant,'' which is derived from the
second construction of subresultant matrix, and show its equivalence
to the recursive subresultant.  In Sect.~\ref{sec:rednessubres}, we
define the ``reduced nested subresultant,'' whose matrix is derived
from the nested subresultant, and show that it is a reduced expression
of the recursive subresultant.  In Sect.~\ref{sec:disc}, we briefly
discuss usage of the reduced nested subresultant in approximate
algebraic computation.

\section{Recursive Polynomial Remainder Sequence (PRS)}
\label{sec:recprs}

First, we review the PRS, then define the recursive PRS.  In the end
of this section, we show a recursive Sturm sequence as an example of
recursive PRS.  We follow definitions and notations by von zur Gathen
and L\"ucking \cite{vzg-luc2003}.  Throughout this paper, let $R$ be
an integral domain and $K$ be its quotient field.  We define a
polynomial remainder sequence as follows.

\begin{defn}[Polynomial Remainder Sequence (PRS)]
  \label{def:prs}
  Let $F$ and \linebreak $G$ be polynomials in $R[x]$ of degree $m$
  and $n$ ($m>n$), respectively.  A sequence
  \begin{equation*}
    (P_1,\ldots,P_l)
  \end{equation*}
  of nonzero polynomials is called a \emph{polynomial remainder
    sequence (PRS)} for $F$ and $G$, abbreviated to $\prs(F,G)$, if it
  satisfies
  \begin{equation*}
    P_1=F,\quad P_2=G,\quad
    \alpha_i P_{i-2} = q_{i-1} P_{i-1} + \beta_i P_{i},
  \end{equation*}
  for $i=3,\ldots,l$, where $\alpha_3,\ldots,\alpha_l,$
  $\beta_3,\ldots,\beta_l$ are elements of $R$ and
  $\deg(P_{i-1})>\deg(P_{i})$.  A sequence
  $((\alpha_3,\beta_3),\ldots,$ $(\alpha_l,\beta_l))$ is called a
  \emph{division rule} for \linebreak $\prs(F,G)$.  If $P_l$ is a
  constant, then the PRS is called \emph{complete}.  \qed
\end{defn}

If $F$ and $G$ are coprime, the last element in the complete PRS for
$F$ and $G$ is a constant.  Otherwise, it equals the GCD of $F$ and
$G$ up to a constant: we have
$\prs(F,G)=(P_{1}=F,P_{2}=G,\ldots,P_{l}=\gamma\cdot \gcd(F,G))$ for
some $\gamma\in R$.  Then, we can calculate new PRS, $\prs(P_{l},
\frac{d}{dx}P_{l})$, and if this PRS ends with a non-constant
polynomial, then calculate another PRS for the last element, and so
on.  By repeating this calculation, we can calculate several PRSs
``recursively'' such that the last polynomial in the last sequence is
a constant.  Thus, we define ``recursive PRS'' as follows.
\begin{defn}[Recursive PRS]
  \label{def:recprs}
  Let $F$ and $G$ be the same as in Definition~\ref{def:prs}. Then, a
  sequence
  \begin{equation*}
    (P_1^{(1)},\ldots,P_{l_1}^{(1)},
    P_1^{(2)},\ldots,P_{l_2}^{(2)},
    \ldots,
    P_1^{(t)},\ldots,P_{l_t}^{(t)})
  \end{equation*}
  of nonzero polynomials is called a \emph{recursive polynomial
    remainder sequence} (recursive PRS) for $F$ and $G$, abbreviated
  to $\rprs(F,G)$, if it satisfies
  \begin{equation*}
    \begin{split}
      & P_1^{(1)} = F,\quad P_2^{(1)}=G,\quad
      P_{l_1}^{(1)}=\gamma_1\cdot\gcd(P_1^{(1)},P_2^{(1)})\;
      \mbox{with $\gamma_1\in R$}, \\
      & (P_1^{(1)},P_2^{(1)},\ldots,P_{l_1}^{(1)})=\prs(P_1^{(1)},P_2^{(1)}),\\
      & P_1^{(k)}=P_{l_{k-1}}^{(k-1)},\quad
      P_2^{(k)}=\frac{d}{dx}P_{l_{k-1}}^{(k-1)},\quad
      P_{l_k}^{(k)}=\gamma_k\cdot\gcd(P_1^{(k)},P_2^{(k)})\;
      \mbox{with $\gamma_k\in R$}, \\
      & (P_1^{(k)},P_2^{(k)},\ldots,P_{l_k}^{(k)})=
      \prs(P_1^{(k)},P_2^{(k)}),
    \end{split}
  \end{equation*}
  for $k=2,\ldots,t$.  If $\alpha_i^{(k)}$, $\beta_i^{(k)}\in R$
  satisfy
  \begin{equation*}
    \alpha_i^{(k)} P_{i-2}^{(k)}
    =
    q_{i-1}^{(k)} P_{i-1}^{(k)}
    + 
    \beta_i^{(k)} P_i^{(k)}
  \end{equation*}
  for $k=1,\ldots,t$ and $i=3,\ldots,l_k$, then a sequence
  $((\alpha_3^{(1)},\beta_3^{(1)}),\ldots,
  (\alpha_{l_t}^{(t)},\beta_{l_t}^{(t)}))$ is called a
  \emph{division rule} for $\rprs(F,G)$. 
  Furthermore, if $P_{l_t}^{(t)}$ is a constant, then the recursive
  PRS is called complete.
  \qed
\end{defn}

\begin{rem}
  \label{re:recprs}
  In this paper, we use the following notations unless otherwise
  defined: for $k=1,\ldots,t$ and $i=1,\ldots,l_k$, let
  $c_i^{(k)}=\lc(P_i^{(k)})$, $n_i^{(k)}=\deg(P_i^{(k)})$ (letters
  $\lc$ and $\deg$ denote the leading coefficient and the degree of
  the polynomial, respectively), $j_0=m$ and $j_k=n_{l_k}^{(k)}$, and
  let $d_i^{(k)}=n_i^{(k)}-n_{i+1}^{(k)}$ for $k=1,\ldots,t$ and
  $i=1,\ldots,l_k-1$.  Furthermore, we represent $P_i^{(k)}(x)$ as
  \[
  P_i^{(k)}(x)=a_{i,n_i^{(k)}}^{(k)}x^{n_i^{(k)}}+\cdots+a_{i,0}^{(k)}x^0,
  \]
  and its ``coefficient vector'' as
  \[
  \bm{p}_i^{(k)}={}^t(a_{i,n_i^{(k)}}^{(k)},\ldots,a_{i,0}^{(k)}).
  \qed
  \]
\end{rem}

As an example of recursive PRS, we calculate Sturm sequences
recursively for calculating the number of real zeros of univariate
polynomial including multiplicities (see Bochnak, Coste and Roy
\cite{boc-cos-roy1998}), as follows.

\begin{exmp}[Recursive Sturm Sequence]
  \label{ex:recsturmseq}
  Let $P(x)=(x+2)^2\{(x-3)(x+1)\}^3$, and calculate the recursive
  Sturm sequence of $P(x)$ as
  \begin{equation*}
    \mathrm{(complete)}\; \rprs(P(x),\frac{d}{dx}P(x)),
  \end{equation*}
  with division rule given by
  \begin{equation*}
    (\alpha_i^{(k)},\beta_i^{(k)})=(1,-1),
  \end{equation*}
  for $k=1,\ldots,t$ and $i=3,\ldots,l_k$.   

  The first sequence 
  $L_1=(P_1^{(1)},\ldots,P_4^{(1)})$ has the following elements:
  \begin{equation*}
    \begin{split}
      P_1^{(1)} &= P(x)=(x+2)^2\{(x-3)(x+1)\}^3,\\
      P_2^{(1)} &=
      \frac{d}{dx}P(x)=8x^7-14x^6-102x^5+80x^4+460x^3+66x^2-558x-324,\\
      P_3^{(1)} &=
      \frac{75}{16}x^6-\frac{45}{16}x^5-60x^4-\frac{225}{8}x^3
      +\frac{3315}{16}x^2
      +\frac{4815}{16}x+\frac{945}{8},\\
      P_4^{(1)} &= \frac{128}{25}x^5-\frac{256}{25}x^4-\frac{256}{5}x^3
      +\frac{1024}{25}x^2+\frac{4224}{25}x+\frac{2304}{25}.
    \end{split}
  \end{equation*}
  The second sequence $L_2=(P_1^{(2)},\ldots,P_4^{(2)})$ has the following
  elements: 
  \begin{equation*}
    \begin{split}
      P_1^{(2)} &= P_4^{(1)} =
      \frac{128}{25}x^5-\frac{256}{25}x^4-\frac{256}{5}x^3
      +\frac{1024}{25}x^2
      +\frac{4224}{25}x+\frac{2304}{25},\\
      P_2^{(2)} &= \frac{d}{dx}P_4^{(1)} =
      \frac{128}{5}x^4-\frac{1024}{25}x^3-\frac{768}{5}x^2
      +\frac{2048}{25}x+\frac{4224}{25},\\
      P_3^{(2)} &= \frac{14848}{625}x^3-\frac{1536}{125}x^2
      -\frac{88576}{625}x-\frac{66048}{625},\\
      P_4^{(2)} &= \frac{12800}{841}x^2-\frac{25600}{841}x
      -\frac{38400}{841}.
    \end{split}
  \end{equation*}
  The last sequence
  $L_3=(P_1^{(3)},\ldots,P_3^{(3)})$ has the following elements:
  \begin{equation*}
    \begin{split}
      P_1^{(3)} &= P_4^{(2)} =
      \frac{12800}{841}x^2-\frac{25600}{841}x -\frac{38400}{841}, \\
      P_2^{(3)} &= \frac{d}{dx}P_4^{(2)} =
      \frac{25600}{841}x-\frac{25600}{841},\\
      P_3^{(3)} &= \frac{51200}{841}.
    \end{split}
  \end{equation*}
  For PRS $L_k$, $k=1,2,3$, define sequences of nonzero real numbers
  $\lambda(L_k,-\infty)$ and $\lambda(L_k,+\infty)$ as
  \begin{equation*}
    \begin{split}
      \lambda(L_k,-\infty) &= 
      ((-1)^{n_1^{(k)}}\lc(P_1^{(k)}),\ldots,
      (-1)^{n_{l_k}^{(k)}}\lc(P_{l_k}^{(k)})),\\ 
      \lambda(L_k,+\infty) &=
      (\lc(P_1^{(k)}),\ldots,
      \lc(P_{l_k}^{(k)})),
    \end{split}
  \end{equation*}
  where $n_i^{(k)}=\deg(P_i^{(k)})$ denotes the degree of $P_i^{(k)}$
  and $\lc(P_i^{(k)})$ denotes the leading coefficients of
  $P_i^{(k)}$.  Then, $\lambda(L_k,-\infty)$ and
  $\lambda(L_k,+\infty)$ for $k=1,2,3$ are
  \begin{equation*}
    \begin{split}
      \lambda(L_1,\pm\infty) &= (1,\pm8,\frac{75}{16},\pm\frac{128}{25}),\\
      \lambda(L_2,\pm\infty) &=
      (\pm\frac{128}{25},\frac{128}{5},\pm\frac{18848}{625},
      \frac{12800}{841}),\\  
      \lambda(L_3,\pm\infty) &=
      (\frac{12800}{841},\pm\frac{25600}{841},\frac{51200}{841}).
    \end{split}
  \end{equation*}
  For a sequence of nonzero real numbers $L=(a_1,\ldots,a_m)$, let
  $V(L)$ be the number of sign variations of the elements of $L$.
  Then, we calculate the number of the real zeros of $P(x)$, including
  multiplicity, as
  \begin{equation*}
    \sum_{k=1}^3\{V(\lambda(L_k,-\infty))-V(\lambda(L_k,+\infty))\}
    = 3+3+2=8.
    \qed
  \end{equation*}
\end{exmp}

\section{Subresultants for Recursive PRS}
\label{sec:subresrecprs}

To make this paper self-contained and to use notations in our
definitions, we first review the fundamental theorem of subresultants,
then discuss subresultants for recursive PRS.

Although the theory of subresultants is established for polynomials
over an integral domain, in what follows, we handle polynomials over a
field for the sake of simplicity.  Let $F$ and $G$ be polynomials in
$K[x]$ such that
\begin{equation}
  \label{eq:fg}
  F(x) = f_m x^m + \cdots + f_0 x^0,\quad
  G(x) = g_n x^n + \cdots + g_0 x^0,\quad
\end{equation}
with $m\ge n>0$.  For a square matrix $M$, we denote its determinant
by $|M|$.

\subsection{Fundamental Theorem of Subresultants}

\begin{defn}[Sylvester Matrix]
  Let $F$ and $G$ be as in \textup{(\ref{eq:fg})}.  The \linebreak
  \emph{Sylvester matrix} of $F$ and $G$, denoted by $N(F,G)$, is an
  $(m+n)\times(m+n)$ matrix constructed from the coefficients of $F$
  and $G$, such that
  \begin{equation*}
    \begin{split}
      N(F,G) &=
      \begin{pmatrix}
        f_m    &        &        & g_n    &        &  \\
        \vdots & \ddots &        & \vdots & \ddots &  \\
        f_0    &        & f_m    & g_0    &        & g_n \\
               & \ddots & \vdots &        & \ddots & \vdots \\
               &        & f_0    &        &        & g_0
      \end{pmatrix}.
      \\
      &\qquad \underbrace{\hspace{1.7cm}}_{n}
      \;\; \underbrace{\hspace{1.5cm}}_{m}
    \end{split}
  \qed
  \end{equation*}
\end{defn}

\begin{defn}[Subresultant Matrix]
  \label{def:subresmat}
  Let $F$ and $G$ be defined as in \textup{(\ref{eq:fg})}.  For $j<n$,
  the \emph{$j$-th subresultant matrix} of $F$ and $G$, denoted by
  $N^{(j)}(F,G)$, is an $(m+n-j)\times(m+n-2j)$ sub-matrix of
  $N(F,G)$ obtained by taking the left $n-j$ columns of
  coefficients of $F$ and the left $m-j$ columns of coefficients of
  $G$, such that
  \begin{equation*}
    \begin{split}
      N^{(j)}(F,G) &=
      \begin{pmatrix}
        f_m    &        &        & g_n    &        &  \\
        \vdots & \ddots &        & \vdots & \ddots &  \\
        f_0    &        & f_m    & g_0    &        & g_n \\
               & \ddots & \vdots &        & \ddots & \vdots \\
               &        & f_0    &        &        & g_0
       \end{pmatrix}.
      \\
      &\qquad \underbrace{\hspace{1.7cm}}_{n-j}
      \;\; \underbrace{\hspace{1.5cm}}_{m-j}
    \end{split}
  \end{equation*}
  Furthermore, define $N_U^{(j)}(F,G)$ as a sub-matrix of
  $N^{(j)}(F,G)$ by deleting the bottom $j+1$ rows.
  \qed
\end{defn}
\begin{defn}[Subresultant]
  Let $F$ and $G$ be defined as in \textup{(\ref{eq:fg})}.  For $j<n$
  and $k=0,\ldots,j$, let $N_k^{(j)}=N_k^{(j)}(F,G)$ (distinguish it
  from $N_U^{(j)}(F,G)$ in the above) be a sub-matrix of
  $N^{(j)}(F,G)$ obtained by taking the top $m+n-2j-1$ rows and the
  $(m+n-j-k)$-th row (note that $N^{(j)}_k(F,G)$ is a square matrix).
  Then, the polynomial
  \begin{equation*}
    \subres_j(F,G)
    =|N^{(j)}_j|x^j+\cdots+|N^{(j)}_0|x^0
  \end{equation*}
  is called the \emph{$j$-th subresultant} of $F$ and $G$. \qed
\end{defn}
\begin{thm}[Fundamental Theorem of Subresultants \cite{bro-tra71}]
  \label{th:fundsubres}
  Let $F$ \linebreak and $G$ be defined as in \textup{(\ref{eq:fg})},
  $(P_1,\ldots,P_k)=\prs(F,G)$ be complete PRS, and
  $((\alpha_3,\beta_3),\ldots,(\alpha_k,\beta_k))$ be its division
  rule.  Let $n_i=\deg(P_i)$ and $c_i=\lc(P_i)$ for $i=1,\ldots,k$,
  and $d_i=n_i-n_{i+1}$ for $i=1,\ldots,k-1$.  Then, we have
  \begin{align}
    \label{eq:subresthm1}
    \subres_j(F,G) &= 0 \quad\textrm{for $0\le j<n_k$},
    \\
    \subres_{n_i}(F,G) &= P_ic_i^{d_{i-1}-1} \nonumber
    \\
    &\qquad\times
    \prod_{l=3}^i
    \left\{
    \left(
      \frac{\beta_l}{\alpha_l}
    \right)^{n_{l-1}-n_i} 
    c_{l-1}^{d_{l-2}+d_{l-1}}
    (-1)^{(n_{l-2}-n_i)(n_{l-1}-n_i)}
    \right\}
    ,
    \\
    \subres_j(F,G) &= 0 \quad\textrm{for $n_i<j<n_{i-1}-1$},
    \\
    \label{eq:subresthm4}
    \subres_{n_{i-1}-1}(F,G) &= P_ic_{i-1}^{1-d_{i-1}} 
    \prod_{l=3}^i
    \Biggl\{\!
    \left(
      \frac{\beta_l}{\alpha_l}
    \right)^{n_{l-1}-n_{i-1}+1}\!\!\! \nonumber
    \\
    &\qquad\times
    c_{l-1}^{d_{l-2}+d_{l-1}}
    (-1)^{(n_{l-2}-n_{i-1}+1)(n_{l-1}-n_{i-1}+1)}\!
    \Biggr\}
    ,
  \end{align}
  for $i=3,\ldots,k$. \qed 
\end{thm}

By the Fundamental Theorem of subresultants, we can express
coefficients of PRS by determinants of matrices whose elements are the
coefficients of initial polynomials.

\subsection{Recursive Subresultants}

We construct ``recursive subresultant matrix'' whose determinants
represent elements of recursive PRS by the coefficients of initial
polynomials.  To help the readers, we first show an example of
recursive subresultant matrix for the recursive Sturm sequence in
Example~\ref{ex:recsturmseq}.
\begin{exmp}[Recursive Subresultant Matrix]
  \label{ex:recsubresmat}
  We express $P(x)$ and \linebreak $\frac{d}{dx}P(x)$ in
  Example~\ref{ex:recsturmseq} by
  \begin{equation*}
    P(x) = f_8 x^8 + \cdots + f_0 x^0,\quad
    \frac{d}{dx}P(x) = g_7 x^7 + \cdots + g_0 x^0.
  \end{equation*}
  Let $\bar{N}^{(1,5)}(F,G)=N^{(5)}(F,G)$, then the matrices
  $\bar{N}_U^{(1,5)}(F,G)$, $\bar{N}_L^{(1,5)}(F,G)$ and\linebreak
  $\bar{N}_L^{'(1,5)}(F,G)$ are given as
  \begin{equation*}
    \begin{split}
      \bar{N}^{(1,5)}(F,G)
      =
      \begin{pmatrix}
        \bar{N}_U^{(1,5)} \\
        \hline
        \bar{N}_L^{(1,5)}
      \end{pmatrix}
      &=
      \begin{pmatrix}
        f_8  &      & g_7 &     &           \\
        f_7  & f_8  & g_6 & g_7 &           \\
        f_6  & f_7  & g_5 & g_6 & g_7 \\
        f_5  & f_6  & g_4 & g_5 & g_6 \\
        \hline
        f_4  & f_5  & g_3 & g_4 & g_5 \\
        f_3  & f_4  & g_2 & g_3 & g_4 \\
        f_2  & f_3  & g_1 & g_2 & g_3 \\
        f_1  & f_2  & g_0 & g_1 & g_2 \\
        f_0  & f_1  &     & g_0 & g_1 \\
        & f_0  &     &     & g_0 \\
      \end{pmatrix}
      ,
      \\
      \bar{N}_L^{'(1,5)}(F,G)
      &
      =
      \begin{pmatrix}
        5f_4 & 5f_5 & 5g_3 & 5g_4 & 5g_5 \\
        4f_3 & 4f_4 & 4g_2 & 4g_3 & 4g_4 \\
        3f_2 & 3f_3 & 3g_1 & 3g_2 & 3g_3 \\
        2f_1 & 2f_2 & 2g_0 & 2g_1 & 2g_2 \\
        f_0  & f_1  &     & g_0 & g_1 \\
      \end{pmatrix}
      ,
    \end{split}
  \end{equation*}
  where horizontal lines in matrices divide them into the upper and
  the lower components.  Note that the matrix $\bar{N}^{'(1,5)}(F,G)$
  is derived from $\bar{N}_L^{(1,5)}(F,G)$ by multiplying the $l$-th
  row by $6-l$ for $l=1,\ldots,5$ and deleting the bottom row.
  Then, the $(2,3)$-th recursive subresultant matrix
  $\bar{N}^{(2,3)}(F,G)$ is constructed as
  \begin{equation}
    \label{eq:recsubresmatex}
    \bar{N}^{(2,3)}(F,G)=
    \left(
      \begin{array}{c|c|c}
        \bar{N}_U^{(1,5)} &  &\\
        \cline{1-2}
                    & \bar{N}_U^{(1,5)} & \\
                    \cline{2-3}
                    &             & \bar{N}_U^{(1,5)}\\
                    \hline
                    &             &0 \cdots 0 \\
                    \cline{3-3}
        \bar{N}_L^{(1,5)} & \bar{N}_L^{'(1,5)} & \\
                    &             & \bar{N}_L^{'(1,5)}\\
                    \cline{2-2}
                    & 0\cdots 0   &\\
      \end{array}
    \right).
  \qed
  \end{equation}
\end{exmp}

\begin{defn}[Recursive Subresultant Matrix]
  \label{def:recsubresmat}
  Let $F$ and $G$ be \linebreak defined as in \textup{(\ref{eq:fg})},
  and let
  $(P_1^{(1)},\ldots,P_{l_1}^{(1)},\ldots,P_1^{(t)},\ldots,P_{l_t}^{(t)})$
  be complete recursive PRS for $F$ and $G$ as in
  Definition~\ref{def:recprs}.  Then, for each pair of numbers $(k,j)$
  with $k=1,\ldots,t$ and $j=j_{k-1}-2,\ldots,0$, define matrix
  $\bar{N}^{(k,j)}=\bar{N}^{(k,j)}(F,G)$ recursively as follows.
  \begin{enumerate}
  \item For $k=1$, let $\bar{N}^{(1,j)}(F,G)=N^{(j)}(F,G)$.
  \item For $k>1$, let $\bar{N}^{(k,j)}(F,G)$ consist of the upper
    block and the lower block, defined as follows:
    \begin{enumerate}
    \item The upper block is partitioned into $(2j_{k-1}-2j-1)\times
      (2j_{k-1}-2j-1)$ blocks with the diagonal blocks filled with
      $\bar{N}_U^{(k-1,j_{k-1})}$, where $\bar{N}_U^{(k-1,j_{k-1})}$
      is a sub-matrix of $\bar{N}^{(k-1,j_{k-1})}(F,G)$ obtained by
      deleting the bottom $j_{k-1}+1$ rows.
    \item Let $\bar{N}_L^{(k-1,j_{k-1})}$ be a sub-matrix of
      $\bar{N}^{(k-1,j_{k-1})}$ obtained by taking the bottom
      $j_{k-1}+1$ rows, and let $\bar{N}_L^{'(k-1,j_{k-1})}$ be a
      sub-matrix of $\bar{N}_L^{(k-1,j_{k-1})}$ by multiplying the
      $(j_{k-1}+1-\tau)$-th rows by $\tau$ for
      $\tau=j_{k-1},\ldots,1$, then by deleting the bottom row.  Then,
      the lower block consists of $j_{k-1}-j-1$ blocks of
      $\bar{N}_L^{(k-1,j_{k-1})}$ such that the leftmost block is
      placed at the top row of the container block and the right-side
      block is placed down by 1 row from the left-side block, then
      followed by $j_{k-1}-j$ blocks of $\bar{N}_L^{'(k-1,j_{k-1})}$
      placed by the same manner as $\bar{N}_L^{(k-1,j_{k-1})}$.
    \end{enumerate}
  \end{enumerate}
As a result, $\bar{N}^{(k,j)}(F,G)$ becomes as shown in
Fig.~\ref{fig:recsubresmat}.  Then, $\bar{N}^{(k,j)}(F,G)$ is called
the \emph{$(k,j)$-th recursive subresultant matrix} of $F$ and $G$.
\qed
\end{defn}

\begin{prop}
  \label{prop:recsubresmat}
  The numbers of rows and columns of $\bar{N}^{(k,j)}(F,G)$, the
  \linebreak $(k,j)$-th recursive subresultant matrix of $F$ and $G$,
  are as follows: for $k=1$ and $j<n$, they are equal to
  \begin{equation}
    \label{eq:recsubresmatorder1}
    m+n-j\quad \textrm{and}\quad m+n-2j,
  \end{equation}
  respectively, and, for $(k,j)=(1,j_1)$ and $k=2,\ldots,t$ and
  $j<j_{k-1}-1$, they are equal to
  \begin{equation}
    \label{eq:recsubresmatrow}
    (m+n-2j_1)
    \left\{
      \prod_{l=2}^{k-1}(2j_{l-1}-2j_l-1)
    \right\}
    (2j_{k-1}-2j-1)
    +j
  \end{equation}
  and
  \begin{equation}
    \label{eq:recsubresmatcol}
    (m+n-2j_1)
    \left\{
      \prod_{l=2}^{k-1}(2j_{l-1}-2j_l-1)
    \right\}
    (2j_{k-1}-2j-1)
    ,
  \end{equation}
  respectively, with $j_0=j_1+1$ for $(k,j)=(1,j_1)$.
\end{prop}
\begin{pf}
  By induction on $k$.  For $k=1$,
  \textup{(\ref{eq:recsubresmatorder1})} immediately follows from Case
  1 of Definition~\ref{def:recsubresmat}, and we also have
  \textup{(\ref{eq:recsubresmatrow})} and
  \textup{(\ref{eq:recsubresmatcol})} for $(k,j)=(1,j_1)$.  Let us
  assume that we have \textup{(\ref{eq:recsubresmatrow})} and
  \textup{(\ref{eq:recsubresmatcol})} for $1,\ldots,k-1$.  Then, we
  calculate the numbers of the rows and columns of
  $\bar{N}^{(k,j)}(F,G)$ as follows.
\begin{enumerate}
\item The numbers of rows of $\bar{N}_L^{(k-1,j_{k-1})}$ and
  $\bar{N}_L^{'(k-1,j_{k-1})}$ are equal to $j_{k-1}+1$ and $j_{k-1}$,
  respectively, thus the number of rows a block which consists of
  $\bar{N}_L^{(k-1,j_{k-1})}$ and $\bar{N}_L^{'(k-1,j_{k-1})}$ in
  $\bar{N}^{(k,j)}(F,G)$ equals
  \begin{equation}
    \label{eq:recsubresmatrow-low}
    2j_{k-1}-j-1.
  \end{equation}
  On the other hand, the number of rows of $\bar{N}_U^{(k-1,j_{k-1})}$
  is equal to $(m+n-2j_1)\{\prod_{l=2}^{k-1}(2j_{l-1}-2j_l-1)\}-1$,
  thus the number of rows of diagonal blocks in
  $\bar{N}^{(k,j)}(F,G)$ is equal to
  \begin{equation}
    \label{eq:recsubresmatrow-up}
    \left\{
      (m+n-2j_1)\prod_{l=2}^{k-1}(2j_{l-1}-2j_l-1)-1
    \right\}
    (2j_{k-1}-2j-1).
  \end{equation}
  By adding (\ref{eq:recsubresmatrow-low}) and
  (\ref{eq:recsubresmatrow-up}), we obtain (\ref{eq:recsubresmatrow}).
\item The number of columns of $\bar{N}^{(k-1,j_{k-1})}(F,G)$ is equal
  to $(m+n-2j_1)$\linebreak $\times\{\prod_{l=2}^{k-1}(2j_{l-1}-2j_l-1)\}$,
  hence the number of columns of $\bar{N}^{(k,j)}(F,G)$ is equal to
  (\ref{eq:recsubresmatcol}).
\end{enumerate}
This proves the proposition.
\qed
\end{pf}

Now, we define recursive subresultants.

\begin{defn}[Recursive Subresultant]
  \label{def:recsubres}
  Let $F$ and $G$ be defined as in \textup{(\ref{eq:fg})}, and let
  $(P_1^{(1)},\ldots,$
  $P_{l_1}^{(1)},\ldots,P_1^{(t)},\ldots,P_{l_t}^{(t)})$ be complete
  recursive PRS for $F$ and $G$ as in Definition~\ref{def:recprs}.
  For $j=j_{k-1}-2,\ldots,0$ and $\tau=j,\ldots,0$, let
  $\bar{N}_\tau^{(k,j)}=\bar{N}_\tau^{(k,j)}(F,G)$ be a sub-matrix of
  the $(k,j)$-th recursive subresultant matrix $\bar{N}^{(k,j)}(F,G)$
  obtained by taking the top
  $(m+n-2j_1)\{\prod_{l=2}^{k-1}(2j_{l-1}-2j_l-1)\}(2j_{k-1}-2j-1)-1$
  rows and the
  $\bigl((m+n-2j_1)\{\prod_{l=2}^{k-1}(2j_{l-1}-2j_l-1)\}(2j_{k-1}-2j-1)
  +j-\tau\bigr)$-th row (note that $\bar{N}_\tau^{(k,j)}$ is a square
  matrix).  Then, the polynomial
  \begin{equation*}
    \recsubres_{k,j}(F,G)
    =|\bar{N}_j^{(k,j)}|x^j+\cdots+|\bar{N}_0^{(k,j)}|x^0
  \end{equation*}
  is called the \emph{$(k,j)$-th recursive subresultant} of
  $F$ and $G$. \qed
\end{defn}

We show the relationship between recursive subresultants and
coefficients in the recursive PRS.

\begin{lem}
  \label{lem:recsubres}
  Let $F$ and $G$ be defined as in \textup{(\ref{eq:fg})}, and let
  $(P_1^{(1)},\ldots,P_{l_1}^{(1)},\ldots,$ 
  $P_1^{(t)},\ldots,P_{l_t}^{(t)})$ 
  be complete recursive PRS for $F$ and $G$ as in
  Definition~\ref{def:recprs}.  For $k=1,\ldots,t-1$, define
  \begin{equation*}
    \begin{split}
      B_k &= (c_{l_k}^{(k)})^{d_{l_k-1}^{(k)}-1}
      \prod_{l=3}^{l_k}
      \Biggl\{
      \left(
        \frac{\beta_l^{(k)}}{\alpha_l^{(k)}}
      \right)^{n_{l-1}^{(k)}-n_{l_k}^{(k)}}
      \!\!
      (c_{l-1}^{(k)})^{(d_{l-2}^{(k)}+d_{l-1}^{(k)})}
      (-1)^{
        (n_{l-2}^{(k)}-n_{l_k}^{(k)})
        (n_{l-1}^{(k)}-n_{l_k}^{(k)})
      }
      \Biggr\}.
    \end{split}
  \end{equation*}
  For $k=2,\ldots,t$ and $j=j_{k-1}-2,\ldots,0$, define
  \begin{equation*}
    \begin{split}
      u_{k,j} &= (m+n-2j_1)
      \left\{
        \prod_{l=2}^{k-1}(2j_{l-1}-2j_l-1)
      \right\}
      (2j_{k-1}-2j-1)
      \\
      &\qquad\mbox{with $u_k=u_{k,j_k}$ and $u_1=m+n-2j_1$}
      ,
      \\
      b_{k,j} &= 2j_{k-1}-2j-1\; \mbox{with $b_k=b_{k,j_k}$ and
      $b_{1,j}=1$ for $j<n$},
      \\
      r_{k,j} &= (-1)^{(u_{k-1}-1)(1+2+\cdots+(b_{k,j}-1))}\;
      \mbox{with $r_k=r_{k,j_k}$ and $r_{1,j}=1$ for $j<n$},
      \\
      \bar{R}_k &= (\bar{R}_{k-1})^{b_k}r_kB_k\; \mbox{with $\bar{R}_0=1$}.
    \end{split}
  \end{equation*}
  Then, for the $(k,j)$-th recursive subresultant of
  $F$ and $G$ with $k=1,\ldots,t$ and $j=j_{k-1}-2,\ldots,0$, we have
  \begin{equation}
    \label{eq:recsubreslem}
      \recsubres_{k,j}(F,G)
      =(\bar{R}_{k-1})^{b_{k,j}}r_{k,j}\cdot\subres_j(P_1^{(k)},P_2^{(k)}).
  \end{equation}
\end{lem}
To prove Lemma~\ref{lem:recsubres}, we prove the following lemma.
\begin{lem}
  \label{lem:recsubresmat}
  For $k=1,\ldots,t$, $j=j_{k-1}-2,\ldots,0$ and $\tau=j,\ldots,0$, we
  have
  \begin{equation*}
      |\bar{N}_\tau^{(k,j)}(F,G)|
      =(\bar{R}_{k-1})^{b_{k,j}}r_{k,j}\,|N_\tau^{(j)}(P_1^{(k)},P_2^{(k)})|.
  \end{equation*}
\end{lem}
\begin{pf}
  By induction on $k$.  For $k=1$, it is obvious from Case 1 in
  Definition~\ref{def:recsubresmat}.  Let us assume that the lemma is valid
  for $1,\ldots,k-1$, then we prove the claim for $k$ by the following
  steps.
\begin{lem}
  \label{lem:recsubresmatm}
  Assume that we have Lemma~\ref{lem:recsubresmat} for $1,\ldots,k-1$.
  Then, for $k$, $j=j_{k-1}-2,\ldots,0$ and $\tau=j,\ldots,0$,
  $\bar{N}^{(k,j)}(F,G)$ can be transformed by certain eliminations
  and permutations on its columns into $M^{(k,j)}(F,G)$ as shown in
  Fig.~\ref{fig:recsubresmatm}, satisfying
  \begin{equation}
    \label{eq:recsubresmatm}
    |\bar{N}^{(k,j)}_\tau(F,G)| =
    ((\bar{R}_{k-2})^{b_{k-1}} r_{k-1})^{b_{k,j}} 
    |M^{(k,j)}_\tau(F,G)|,
  \end{equation}
  where $M_\tau^{(k,j)}(F,G)$ is a sub-matrix of $M^{(k,j)}(F,G)$
  obtained by the same manner as we have obtained
  $\bar{N}_\tau^{(k,j)}(F,G)$ from $\bar{N}^{(k,j)}(F,G)$ in
  Definition~\ref{def:recsubres}.
\end{lem}
\begin{pf}  
  By the induction hypothesis, for $\tau'=j_{k-1},\ldots,0$, we
  have
  \begin{equation*}
    |\bar{N}_{\tau'}^{(k-1,j_{k-1})}(F,G)|=
    (\bar{R}_{k-2})^{b_{k-1}}r_{k-1}\,
    |N_{\tau'}^{(j_{k-1})}(P_1^{(k-1)},P_2^{(k-1)})|. 
  \end{equation*}
  Let $\bar{N}'^{(k,j)}(F,G)$ be a matrix defined as
  \begin{equation*}
    \bar{N}'^{(k,j)}(F,G)=
    \begin{pmatrix}
      \bar{N}_U^{(k,j)} \\
      \hline
      \bar{N}_L'^{(k,j)}
    \end{pmatrix}
    ,
  \end{equation*}
  where $\bar{N}_U^{(k,j)}$ and $\bar{N}_L'^{(k,j)}$ are defined as in
  Definition~\ref{def:recsubresmat}.  Furthermore, let\linebreak
  $N'^{(j_{k-1})}(P_1^{(k-1)},P_2^{(k-1)})$ be defined as
  $N^{(j_{k-1})}(P_1^{(k-1)},P_2^{(k-1)})$ with the
  $(j_{k-2}+1-\tau)$-th row multiplied by $\tau$ for
  $\tau=j_{k-1},\ldots,1$, then by deleting the bottom row,
  $N_U^{(j_{k-1})}(P_1^{(k-1)},P_2^{(k-1)})$ be a sub-matrix of
  $N^{(j_{k-1})}(P_1^{(k-1)},P_2^{(k-1)})$ and\linebreak
  $N'^{(j_{k-1})}(P_1^{(k-1)},P_2^{(k-1)})$ by taking the top
  $j_{k-2}-j_{k-1}$ rows, and\linebreak
  $N_L^{(j_{k-1})}(P_1^{(k-1)},P_2^{(k-1)})$ and
  $N_L'^{(j_{k-1})}(P_1^{(k-1)},P_2^{(k-1)})$ be sub-matrices of\linebreak
  $N^{(j_{k-1})}(P_1^{(k-1)},P_2^{(k-1)})$ and
  $N'^{(j_{k-1})}(P_1^{(k-1)},P_2^{(k-1)})$, respectively, by
  eliminating the top $j_{k-2}-j_{k-1}$ rows.
  
  Then, by certain eliminations and exchanges on columns, we can
  transform $\bar{N}^{(k-1,j_{k-1})}(F,G)$ and
  $\bar{N}'^{(k-1,j_{k-1})}(F,G)$ to
  \begin{equation}
    \label{eq:recsubresmatd}
    \begin{split}
      & D^{(k-1,j_{k-1})}(F,G)
      \\
      =&
      \begin{pmatrix}
        W_{k-1} & \vline & \bm{0} \\
        \hline
        *       & \vline & N^{(j_{k-1})}(P_1^{(k-1)},P_2^{(k-1)}) 
      \end{pmatrix}
      =
      \begin{pmatrix}
        W_{k-1} & \vline & \bm{0} \\
        \hline
        *      & \vline & {N_U}^{(j_{k-1})}(P_1^{(k-1)},P_2^{(k-1)}) \\
        \hline
        *      & \vline & {N_L}^{(j_{k-1})}(P_1^{(k-1)},P_2^{(k-1)}) 
      \end{pmatrix}
      ,
      \\
      & {D'}^{(k-1,j_{k-1})}(F,G)
      \\
      =&
      \begin{pmatrix}
        W_{k-1} & \vline & \bm{0} \\
        \hline
        *       & \vline & {N'}^{(j_{k-1})}(P_1^{(k-1)},P_2^{(k-1)}) 
      \end{pmatrix}
      =
      \begin{pmatrix}
        W_{k-1} & \vline & \bm{0} \\
        \hline
        *      & \vline & {N_U}^{(j_{k-1})}(P_1^{(k-1)},P_2^{(k-1)}) \\
        \hline
        *      & \vline & {N_L}'^{(j_{k-1})}(P_1^{(k-1)},P_2^{(k-1)}) 
      \end{pmatrix}
      ,
    \end{split}
  \end{equation}
  respectively, satisfying
  \begin{equation*}
    \begin{split}
      |W_{k-1}| &= 1,
      \\
      |{D}_{\tau'}^{(k-1,j_{k-1})}(F,G)| &= 
      (\bar{R}_{k-2})^{b_{k-1}}r_{k-1}
      |{N}_{\tau'}^{(j_{k-1})}(P_1^{(k-1)},P_2^{(k-1)})|,
      \\
      |{D'}_{\tau'}^{(k-1,j_{k-1})}(F,G)| &= 
      (\bar{R}_{k-2})^{b_{k-1}}r_{k-1}
      |{N'}_{\tau'}^{(j_{k-1})}(P_1^{(k-1)},P_2^{(k-1)})|,  
    \end{split}
  \end{equation*}
  where ${D}_{\tau'}^{(k-1,j_{k-1})}(F,G)$ and
  ${D'}_{\tau'}^{(k-1,j_{k-1})}(F,G)$ are sub-matrices of\linebreak
  ${D}^{(k-1,j_{k-1})}(F,G)$ and ${D'}^{(k-1,j_{k-1})}(F,G)$,
  respectively, obtained by the same manner as we have obtained
  $\bar{N}_{\tau'}^{(k-1,j_{k-1})}(F,G)$ from
  $\bar{N}^{(k-1,j_{k-1})}(F,G)$ (see Definition~\ref{def:recsubres}).
  
  Therefore, by the above transformations on the columns in each
  column blocks in $\bar{N}^{(k,j)}(F,G)$ as shown in
  Fig.~\ref{fig:recsubresmat}, we obtain $M^{(k,j)}(F,G)$ as shown in
  Fig.~\ref{fig:recsubresmatm}, where
  ${N_U}^{(j_{k-1})}(P_1^{(k-1)},P_2^{(k-1)})$,
  ${N_L}^{(j_{k-1})}(P_1^{(k-1)},P_2^{(k-1)})$ and \linebreak
  ${N'_L}^{(j_{k-1})}(P_1^{(k-1)},P_2^{(k-1)})$ are abbreviated to
  ${N_U}^{(j_{k-1})}$, ${N_L}^{(j_{k-1})}$ and ${N'_L}^{(j_{k-1})}$,
  respectively, satisfying (\ref{eq:recsubresmatm}) because
  $M^{(k,j)}(F,G)$ has $b_{k,j}=2j_{k-1}-2j-1$ column blocks that have
  ${D}^{(k-1,j_{k-1})}(F,G)$ or
  ${D'}^{(k-1,j_{k-1})}(F,G)$.  This proves the lemma. \qed
\end{pf}

\begin{lem}
  \label{lem:recsubresmatmbar}
  For $k$, $j=j_{k-1}-2,\ldots,0$ and $\tau=j,\ldots,0$,
  $M^{(k,j)}(F,G)$ can be transformed by certain eliminations and
  permutations on its columns into $\bar{M}^{(k,j)}(F,G)$ as shown
  in Fig.~\ref{fig:recsubresmatmbar}, satisfying
  \begin{equation}
    \label{eq:recsubresmatmbar}
    |M_\tau^{(k,j)}(F,G)|= (B_{k-1})^{b_{k,j}}
    |\bar{M}_\tau^{(k,j)}(F,G)|,
  \end{equation}
  where $\bar{M}_\tau^{(k,j)}$ is a sub-matrix of $\bar{M}^{(k,j)}$
  obtained by the same manner as we have obtained
  $\bar{N}_{\tau'}^{(k-1,j_{k-1})}$ from $\bar{N}^{(k-1,j_{k-1})}$
  in Definition~\ref{def:recsubres}.
\end{lem}
\begin{pf}
  For ${N'}^{(j_{k-1})}(P_1^{(k-1)},P_2^{(k-1)})$ (defined as in the
  proof of Lemma~\ref{lem:recsubresmatm}) and
  $\tau''=j_{k-1}-1,\ldots,0$, let
  ${N'}^{(j_{k-1})}_{\tau''}(P_1^{(k-1)},P_2^{(k-1)})$ be a sub-matrix
  of \linebreak ${N'}^{(j_{k-1})}(P_1^{(k-1)},P_2^{(k-1)})$ obtained
  by taking the top $2(n_1^{(k-1)}-j_{k-1}-1)$ rows and the
  $(2n_1^{(k-1)}-2-j_{k-1}-\tau'')$-th row.  Then, by the Fundamental
  Theorem of subresultants (Theorem~\ref{th:fundsubres}), we have
  \begin{equation*}
    \begin{split}
      &\quad
      |N^{(j_{k-1})}_{j_{k-1}}(P_1^{(k-1)},P_2^{(k-1)})|x^{j_{k-1}}
      +\cdots+
      |N_0^{(j_{k-1})}(P_1^{(k-1)},P_2^{(k-1)})|x^0
      \\
      &=
      \subres_{j_{k-1}}(P_1^{(k-1)},P_2^{(k-1)})
      =B_{k-1}P_{l_{k-1}}^{(k-1)}
      =B_{k-1}P_1^{(k)},
      \\
      &\quad
      |{N'}^{(j_{k-1})-1}_{j_{k-1}}(P_1^{(k-1)},P_2^{(k-1)})|x^{j_{k-1}-1}
      +\cdots+
      |{N'}_0^{(j_{k-1})}(P_1^{(k-1)},P_2^{(k-1)})|x^0
      \\
      &=
      \frac{d}{dx}
      \subres_{j_{k-1}}(P_1^{(k-1)},P_2^{(k-1)})
      =B_{k-1}\frac{d}{dx}P_{l_{k-1}}^{(k-1)}
      =B_{k-1}P_2^{(k)},
    \end{split}
  \end{equation*}
  hence, for $\tau'=j_{k-1},\ldots,0$ and $\tau''=j_{k-1}-1,\ldots,0$,
  we have 
  \begin{equation*}
    \begin{split}
      |N^{(j_{k-1})}_{\tau'}(P_1^{(k-1)},P_2^{(k-1)})| &=
      B_{k-1}a^{(k)}_{1,\tau'},
      \\
      |{N'}^{(j_{k-1})}_{\tau''}(P_1^{(k-1)},P_2^{(k-1)})| &=
      B_{k-1}a^{(k)}_{2,\tau''},
    \end{split}
  \end{equation*}
  where $a^{(k)}_{i,j}$ represents the coefficient of degree $j$ of
  $P^{(k)}_i$ (see Remark~\ref{re:recprs}).  Therefore, by certain
  eliminations and exchanges on columns, we can transform
  $N^{(j_{k-1})}(P_1^{(k-1)},P_2^{(k-1)})$ and
  ${N'}^{(j_{k-1})}(P_1^{(k-1)},P_2^{(k-1)})$ into
  \begin{equation*}
    \left(
      \begin{array}{c|c}
        \bar{N}_U^{(j_{k-1})} & \bm{0} \\
        \hline
        \;\;* & \bm{p}^{(k)}_1
      \end{array}
    \right)
    \quad\textrm{and}\quad
    \left(
      \begin{array}{c|c}
        \bar{N}_U^{(j_{k-1})} & \bm{0} \\
        \hline
        \;\;* & \bm{p}^{(k)}_2
      \end{array}
    \right)
    ,
  \end{equation*}
  respectively, satisfying $|\bar{N}_U^{(j_{k-1})}|=1$ (see
  Remark~\ref{re:recprs} for the notation of ``coefficient vectors'').
  By these transformations, we can transform
  ${D}^{(k-1,j_{k-1})}(F,G)$ and ${D'}^{(k-1,j_{k-1})}(F,G)$ in
  (\ref{eq:recsubresmatd}) to
  \begin{equation*}
    \begin{split}
      \bar{D}^{(k-1,j_{k-1})}(F,G) &=
      \left(
        \begin{array}{c|c|c}
          W_{k-1} & \multicolumn{2}{c}{\bm{0}}\\
          \cline{1-2}
          & \bar{N}_U^{(j_{k-1})} &       \\
          \cline{2-3}
          \multicolumn{2}{l|}{\;\;*}
          & \bm{p}_1^{(k)}
        \end{array}
      \right),
      \\
      \bar{D'}^{(k-1,j_{k-1})}(F,G) &=
      \left(
        \begin{array}{c|c|c}
          W_{k-1} & \multicolumn{2}{c}{\bm{0}}\\
          \cline{1-2}
          & \bar{N}_U^{(j_{k-1})} &       \\
          \cline{2-3}
          \multicolumn{2}{l|}{\;\;*}
          & \bm{p}_2^{(k)}
        \end{array}
      \right),
    \end{split}
  \end{equation*}
  respectively, satisfying
  \begin{equation*}
    \begin{split}
      |\bar{N}_U^{(j_{k-1})}| &= 1,
      \\
      |\bar{D}_{\tau'}^{(k-1,j_{k-1})}(F,G)| &= B_{k-1}
      |{D}_{\tau'}^{(k-1,j_{k-1})}(F,G)|,
      \\
      |\bar{D}_{\tau'}'^{(k-1,j_{k-1})}(F,G)| &= B_{k-1}
      |{D}_{\tau'}'^{(k-1,j_{k-1})}(F,G)|,
    \end{split}
  \end{equation*}
  where ${D}_{\tau'}^{(k-1,j_{k-1})}$, ${D}_{\tau'}'^{(k-1,j_{k-1})}$,
  $\bar{D}_{\tau'}^{(k-1,j_{k-1})}$ and
  $\bar{D}_{\tau'}'^{(k-1,j_{k-1})}$ are sub-matrices of\linebreak
  ${D}^{(k-1,j_{k-1})}$, ${D}'^{(k-1,j_{k-1})}$,
  $\bar{D}^{(k-1,j_{k-1})}$ and $\bar{D}'^{(k-1,j_{k-1})}$,
  respectively, obtained by the same manner as we have obtained
  $\bar{N}_{\tau'}^{(k-1,j_{k-1})}$ from $\bar{N}^{(k-1,j_{k-1})}$.
  Therefore, by the above eliminations on the columns in each column
  blocks, we can transform $M^{(k,j)}(F,G)$ to $\bar{M}^{(k,j)}(F,G)$
  as shown in Fig.~\ref{fig:recsubresmatmbar} satisfying
  (\ref{eq:recsubresmatmbar}) because $M^{(k,j)}(F,G)$ has
  $b_{k,j}=2j_{k-1}-2j-1$ column blocks that have
  ${D}^{(k-1,j_{k-1})}(F,G)$ or ${D'}^{(k-1,j_{k-1})}(F,G)$. This
  proves the lemma. \qed
\end{pf}
\emph{Proof of Lemma~\ref{lem:recsubresmat} (continued).}  By
exchanges on column blocks, we can transform $\bar{M}^{(k,j)}(F,G)$ to
$\hat{M}^{(k,j)}(F,G)$ as shown in Fig.~\ref{fig:recsubresmathat},
with
  \begin{equation}
    \label{eq:recsubresmatmhat}
    |\bar{M}_\tau^{(k,j)}(F,G)|=
    r_{k,j}\ |\hat{M}_\tau^{(k,j)}(F,G)|,
  \end{equation}
  where $\hat{M}_\tau^{(k,j)}$ is a sub-matrix of $\hat{M}^{(k,j)}$
  obtained by the same manner as we have obtained
  $\bar{N}_{\tau}^{(k,j)}$ from $\bar{N}^{(k,j)}$, because the
  $(u_{k,j}-(l-1)u_{k-1})$-th column in $\bar{M}^{(k,j)}(F,G)$ was
  moved to the $(u_{k,j}-(l-1))$-th column in $\hat{M}^{(k,j)}(F,G)$
  for $l=1,\ldots,b_{k,j}$.  (Note that $\hat{M}^{(k,j)}(F,G)$ is a
  block lower triangular matrix.)  Then, we have
  \begin{equation}
    \label{eq:recsubresmatsubres}
    |\hat{M}_\tau^{(k,j)}(F,G)|
    =
    |N_\tau^{(j)}(P_1^{(k)},P_2^{(k)})|,
  \end{equation}
  because we have $|W_{k-1}|=|\bar{N}_U^{(j_{k-1})}|=1$
  and the lower-right block of
  $\bm{p}_1^{(k)}$s and $\bm{p}_2^{(k)}$s in $\hat{M}^{(k,j)}(F,G)$ is
  equal to $N^{(j)}(P_1^{(k)},P_2^{(k)})$.

  Finally, from (\ref{eq:recsubresmatm}), (\ref{eq:recsubresmatmbar}),
  (\ref{eq:recsubresmatmhat}) and (\ref{eq:recsubresmatsubres}), we
  have
  \[
  \begin{split}
    |\bar{N}_\tau^{(k,j)}(F,G)| &=
    ((\bar{R}_{k-2})^{b_{k-1}} r_{k-1})^{b_{k,j}}
    (B_{k-1})^{b_{k,j}}
    r_{k,j}
    |N_\tau^{(j)}(P_1^{(k)},P_2^{(k)})|
    \\
    &=
    (\bar{R}_{k-1})^{b_{k,j}}r_{k,j}|N_\tau^{(j)}(P_1^{(k)},P_2^{(k)})|,
  \end{split}
  \]
  which proves the lemma. \qed
\end{pf}

\begin{thm}
  With the same conditions as in  Lemma~\ref{lem:recsubres}, and for
  $k=1,\ldots,t$ and $i=3,4,\ldots,l_k$, we have
  \begin{align}
    \label{eq:recsubresthm1}
    &\recsubres_{k,j}(F,G) = 0 \quad\mbox{for $0\le j<n_{l_k}^{(k)}$},
    \\
    &\recsubres_{k,n_i^{(k)}}(F,G) = P_i^{(k)}
    (c_i^{(k)})^{d_{i-1}^{(k)}-1} 
    (R_{k-1})^{b_k,n_i^{(k)}} r_{k,n_i^{(k)}}\nonumber
    \\
    & \qquad \times \prod_{l=3}^{i} \left\{ \left(
        \frac{\beta_l^{(k)}}{\alpha_l^{(k)}}
      \right)^{n_{l-1}^{(k)}-n_i^{(k)}}
      (c_{l-1}^{(k)})^{(d_{l-2}^{(k)}+d_{l-1}^{(k)})} (-1)^{
        (n_{l-2}^{(k)}-n_{i}^{(k)}) (n_{l-1}^{(k)}-n_{i}^{(k)}) }
    \right\},
    \\
    &\recsubres_{k,j}(F,G) = 0 \quad\mbox{for
      $n_i^{(k)}<j<n_{i-1}^{(k)}-1$},
    \\
    \label{eq:recsubresthm4}
    &\recsubres_{k,n_{i-1}^{(k)}-1}(F,G) \nonumber
    \\
    &\qquad
    = P_i^{(k)}
    (c_{i-1}^{(k)})^{1-d_{i-1}^{(k)}}
    (R_{k-1})^{b_k,n_{i-1}^{(k)}-1} r_{k,n_{i-1}^{(k)}-1}
    \cdot \prod_{l=3}^{i}
    \Biggl\{
      \left(
        \frac{\beta_l^{(k)}}{\alpha_l^{(k)}}
      \right)^{n_{l-1}^{(k)}-n_{i-1}^{(k)}+1} \nonumber
      \\
      &\qquad\qquad\qquad\qquad\quad
      \times
      (c_{l-1}^{(k)})^{(d_{l-2}^{(k)}+d_{l-1}^{(k)})} (-1)^{
        (n_{l-2}^{(k)}-n_{i-1}^{(k)}+1)
        (n_{l-1}^{(k)}-n_{i-1}^{(k)}+1) }
    \Biggr\}.
  \end{align}
\end{thm}
\begin{pf}
  By substituting $\subres_j(P_1^{(k)},P_2^{(k)})$ in
  (\ref{eq:recsubreslem}) by
  (\ref{eq:subresthm1})--(\ref{eq:subresthm4}), we obtain
  (\ref{eq:recsubresthm1})--(\ref{eq:recsubresthm4}), respectively.
  \qed
\end{pf}

We show an example of the proof of Lemma~\ref{lem:recsubres} for the
recursive subresultant matrix in Example~\ref{ex:recsubresmat}.
\begin{exmp}
  (Continued from Example~\ref{ex:recsubresmat}.)  Since we have
  $\bar{N}^{(1,5)}(F,G)=N^{(5)}(F,G)$, we can regard
  $\bar{N}^{(1,5)}(F,G)$ and $\bar{N'}^{(1,5)}(F,G)$ as
  $D^{(1,5)}(F,G)$ and ${D'}^{(1,5)}(F,G)$ in
  (\ref{eq:recsubresmatd}), respectively.  Then, by eliminations and
  exchanges of columns as shown in Lemma~\ref{lem:recsubresmatmbar},
  we can transform $\bar{N}^{(1,5)}(F,G)=
  \begin{pmatrix}
    \bar{N}_U^{(1,5)} \\
    \hline
    \bar{N}_L^{(1,5)}
  \end{pmatrix}
  $ and
  $\bar{N}^{'(1,5)}(F,G)=
  \begin{pmatrix}
    \bar{N}_U^{(1,5)} \\
    \hline
    \bar{N}_L^{'(1,5)}
  \end{pmatrix}
  $ in (\ref{eq:recsubresmatex}) to
  $\bar{D}^{(1,5)}(F,G)$ and $\bar{D}^{'(1,5)}(F,G)$, respectively, as
  \begin{equation*}
    \bar{D}^{(1,5)}(F,G)
      =
      \left(
        \begin{array}{c|c}
          \bar{N}_U^{(5)} & \bm{0} \\
          \hline
          *  & \bm{p}_1^{(2)}\\
        \end{array}
      \right)
      ,
      \quad
      \bar{D}^{'(1,5)}(F,G)
      =
      \left(
        \begin{array}{c|c}
          \bar{N}_U^{(5)} & \bm{0} \\
          \hline
          *  & \bm{p}_2^{(2)}\\
        \end{array}
      \right)
      ,
    \end{equation*}
  with $|\bar{N}_U^{(5)}|=1$ and $B_1=(a_{2,7}^{(1)})^2 (a_{3,6}^{(1)})^2$.
  Therefore, by the above transformations
  of columns in each column blocks in
  $\bar{N}^{(2,1)}(F,G)$, we have
  \begin{equation*}
    \bar{M}^{(2,3)}(F,G)=
    \left(
      \begin{array}{c|c|c|c|c|c}
        \bar{N}_U^{(5)} & \bm{0} & & & &\\
        \cline{1-3}
                        &         &  \bar{N}_U^{(5)} & \bm{0} & &\\
        \cline{3-5}
                        &         &                  &         &
        \bar{N}_U^{(5)} & \bm{0} \\
        \hline
                        &         &                  &         & \multicolumn{2}{c}{0\cdots 0} \\
        \cline{5-6}
        *               & \bm{p}_1^{(2)} & * & \bm{p}_2^{(2)} &  &\\
                        &         &                  &          & * & \bm{p}_2^{(2)}\\
        \cline{3-4}
                        &         & \multicolumn{2}{c|}{0\cdots 0}
        & & \\
      \end{array}
    \right),
  \end{equation*}
  satisfying
  $|\bar{N}_\tau^{(2,3)}(F,G)|=(B_1)^3\ |\bar{M}_\tau^{(2,3)}(F,G)|$
  for $\tau=3,\ldots,0$. 
  Furthermore, by exchanges of columns, we can transform
  $\bar{M}^{(2,3)}(F,G)$ to $\hat{M}^{(2,3)}(F,G)$ as
  \begin{equation*}
    \begin{split}
    & \hat{M}^{(2,3)}(F,G)
    \\
    =&
    \left(
      \begin{array}{c|c|c|c|c|c}
        \bar{N}_U^{(5)}   \\
        \cline{1-2}
                        & \bar{N}_U^{(5)} \\
        \cline{2-3}
                        &      & \bar{N}_U^{(5)}    \\
        \hline
                        &         & {0\cdots 0}      &         & & 0  \\
        \cline{3-3}\cline{6-6}
        *               &  * & * & \bm{p}_1^{(2)} & \bm{p}_2^{(2)} &\\
                        &         &                  &          &  & \bm{p}_2^{(2)}\\
        \cline{2-2}\cline{5-5}
                        & {0\cdots 0}     &  &
        & 0 & \\
      \end{array}
    \right)
    =
    \left(
      \begin{array}{c|c|c|c}
        \bar{N}_U^{(5)}   \\
        \cline{1-2}
                        & \bar{N}_U^{(5)} \\
        \cline{2-3}
        \multicolumn{2}{c|}{}  & \bar{N}_U^{(5)}    \\
        \hline
        \multicolumn{3}{c|}{*} & N^{(3)}(P_1^{(2)},P_2^{(2)})
      \end{array}
    \right),
    \end{split}
  \end{equation*}
  satisfying
  $|\bar{M}_\tau^{(2,3)}(F,G)|=r_{2,3}\ |\hat{M}_\tau^{(2,3)}(F,G)|
  =r_{2,3}\ |N_\tau^{(3)}(P_1^{(2)},P_2^{(2)})|$.
  Therefore, we have
  \begin{equation*}
    |\bar{N}_\tau^{(2,3)}(F,G)|=(B_1)^3r_{2,3}\ |N_\tau^{(3)}(P_1^{(2)},P_2^{(2)})|
    = (R_1)^3 r_{2,3}\ |N_\tau^{(3)}(P_1^{(2)},P_2^{(2)})|,
  \end{equation*}
  for $\tau=3,\ldots,0$, and we have
  \begin{equation*}
    \recsubres_{2,3}(F,G)=(R_1)^3 r_{2,3}\cdot\subres_3(P_1^{(2)},P_2^{(2)})=
    \{(a_{2,7}^{(1)})^2(a_{3,6}^{(1)})^2\}^3(a_{2,4}^{(2)})^2\ P_3^{(2)}.
  \qed
  \end{equation*}
\end{exmp}

\section{Nested Subresultants}
\label{sec:nessubres}

As we have seen in the above, the recursive subresultant can represent
the coefficients of the elements in recursive PRS.  However, the size
of the recursive subresultant matrix becomes larger rapidly as the
recursion of the recursive PRS deepens, thus making use
of the recursive subresultant matrix become inefficient.

To overcome this problem, we should introduce other representations
for the subresultant that are equivalent to the recursive
subresultant, and more suitable for efficient computations.  The
nested subresultant matrix is a subresultant matrix whose elements are
again determinants of certain subresultant matrices (or even the
nested subresultant matrices), and the nested subresultant is a
subresultant whose coefficients are determinants of the nested
subresultant matrices.

Note that the nested subresultant is mainly used to show the
relationship between the recursive subresultant and the reduced nested
subresultant that will be defined in the next section.

We show an example of a nested subresultant matrix.

\begin{exmp}
  Let $F(x)$ and $G(x)$ be defined as
  \begin{equation*}
    \begin{split}
      F(x) &= a_6 x^6+a_5 x^5+\cdots+a_0,\quad a_6\ne 0,\\
      G(x) &= b_5 x^5+b_4 x^4+\cdots+b_0,\quad b_5\ne 0.
    \end{split}
  \end{equation*}
  Let $\prs(F,G)=(P_1^{(1)}=F,\ P_2^{(1)}=G,\ P_3^{(1)}=\gcd(F,G))$
  with $\deg(P_3^{(1)})=4$, and let us consider recursive PRS for $F$
  and $G$.
  
  Let $P_1^{(2)}=P_3^{(1)}$, $P_2^{(2)}=\frac{d}{dx}P_3^{(1)}$, and
  calculate a subresultant of degree $1$, which corresponds to
  $P_4^{(2)}$.  By the Fundamental Theorem of subresultants
  (Theorem~\ref{th:fundsubres}), we have
  \begin{equation*}
    \begin{split}
      \subres_4(F,G) &= A_4 x^4 + A_3 x^3 + A_2 x^2 + A_1 x +A_0,
      \\
      \frac{d}{dx}\subres_4(F,G) &= 4 A_4 x^3 + 3 A_3 x^2 + 2 A_2 x +
      A_1,
    \end{split}
  \end{equation*}
  where
  \begin{equation}
    \label{eq:exaj}
    A_j=|N_j^{(4)}(F,G)|
  \end{equation}
  for $j=0,\ldots,4$ with $N_k^{(j)}(F,G)$ as in
  Definition~\ref{def:subresmat}. 

  Then, we can express the subresultant matrix
  $N^{(2)}(\subres_4(F,G),\frac{d}{dx}\subres_4(F,G))$ as
  \begin{equation}
    \label{eq:ajmat}
    N^{(2)}(\subres_4(F,G),\frac{d}{dx}\subres_4(F,G))=
    \begin{pmatrix}
      A_4 & 4A_4 & \\
      A_3 & 3A_3 & 4A_4 \\
      A_2 & 2A_2 & 3A_3 \\
      A_1 & A_1  & 2A_2 \\
      A_0 &      & A_1
    \end{pmatrix}
    ,
  \end{equation}
  and the subresultant
  $\subres_2(\subres_4(F,G),\frac{d}{dx}\subres_4(F,G))$ as
  \begin{multline}
    \label{eq:ajdet}
    \subres_2(\subres_4(F,G),\frac{d}{dx}\subres_4(F,G))
    \\
    =
    \begin{vmatrix}
      A_4 & 4A_4 &      \\
      A_3 & 3A_3 & 4A_4 \\
      A_2 & 2A_2 & 3A_3 \\
    \end{vmatrix}
    x^2 +
    \begin{vmatrix}
      A_4 & 4A_4 &      \\
      A_3 & 3A_3 & 4A_4 \\
      A_1 &  A_1 & 2A_2 \\
    \end{vmatrix}
    x + 
    \begin{vmatrix}
      A_4 & 4A_4 &      \\
      A_3 & 3A_3 & 4A_4 \\
      A_0 &      &  A_1 \\
    \end{vmatrix}
    ,
  \end{multline}
  respectively, with $A_j$ as in (\ref{eq:exaj}).  We see that the
  elements in (\ref{eq:ajmat}) are minors of subresultant matrix,
  hence the coefficients in (\ref{eq:ajdet}) is ``nested'' expression
  of determinants.  \qed
\end{exmp}
\begin{defn}[Nested Subresultant Matrix]
  \label{def:nessubresmat}
  Let $F$ and $G$ be defined as in (\ref{eq:fg}), and let
  $(P_1^{(1)},\ldots,P_{l_1}^{(1)},\ldots,P_1^{(t)},\ldots,P_{l_t}^{(t)})$
  be complete recursive PRS for $F$ and $G$ as in
  Definition~\ref{def:recprs}.  Then, for each pair of numbers
  $(k,j)$ with $k=1,\ldots,t$ and $j=j_{k-1}-2,\ldots,0$, define
  matrix $\tilde{N}^{(k,j)}(F,G)$ recursively as follows.
  \begin{enumerate}
  \item For $k=1$, let $\tilde{N}^{(1,j)}(F,G)=N^{(j)}(F,G)$.
  \item For $k>1$, let
    \begin{equation*}
      \tilde{N}^{(k,j)}(F,G)=
      N^{(j)}
      \left(
        \tilde{\subres}_{k-1,j_{k-1}}(F,G),
        \frac{d}{dx}\tilde{\subres}_{k-1,j_{k-1}}(F,G)
      \right),
    \end{equation*}
    where $\tilde{\subres}_{k-1,j_{k-1}}(F,G)$ is defined by
    Definition~\ref{def:nessubres}.  Then, $\tilde{N}^{(k,j)}(F,G)$ is
    called the \emph{$(k,j)$-th nested subresultant matrix of $F$ and
      $G$}.  \qed
  \end{enumerate}
\end{defn}

\begin{defn}[Nested Subresultant]
  \label{def:nessubres}
  Let $F$ and $G$ be defined as in \textup{(\ref{eq:fg})}, and let
  $(P_1^{(1)},\ldots,$
  $P_{l_1}^{(1)},\ldots,P_1^{(t)},\ldots,P_{l_t}^{(t)})$ be complete
  recursive PRS for $F$ and $G$ as in Definition~\ref{def:recprs}.
  For $j=j_{k-1}-2,\ldots,0$ and $\tau=j,\ldots,0$, let
  $\tilde{N}_\tau^{(k,j)}=\tilde{N}_\tau^{(k,j)}(F,G)$ be a sub-matrix
  of the $(k,j)$-th nested subresultant matrix
  $\tilde{N}^{(k,j)}(F,G)$ obtained by taking the top
  $n_1^{(k)}+n_2^{(k)}-2j-1$ rows and the
  $(n_1^{(k)}+n_2^{(k)}-j-\tau)$-th row (note that
  $\tilde{N}_\tau^{(k,j)}$ is a square matrix).  Then, the polynomial
  \begin{equation*}
    \nessubres_{k,j}(F,G)
    =|\tilde{N}_j^{(k,j)}|x^j+\cdots+|\tilde{N}_0^{(k,j)}|x^0
  \end{equation*}
  is called the \emph{$(k,j)$-th nested subresultant} of
  $F$ and $G$. \qed
\end{defn}

We show the relationship between the nested subresultant and the
recursive subresultant.

\begin{lem}
  \label{lem:nessubres}
  Let $F$ and $G$ be defined as in \textup{(\ref{eq:fg})}, and let
  $(P_1^{(1)},\ldots,P_{l_1}^{(1)},\ldots,$,
  $P_1^{(t)},\ldots,P_{l_t}^{(t)})$ be complete recursive PRS for $F$
  and $G$ as in Definition~\ref{def:recprs}.  For $k=1,\ldots,t-1,$
  define $B_k$, and for $k=2,\ldots,t$ and $j=j_{k-1}-2,\ldots,0$,
  define $b_{k,j}$ as in Lemma~\ref{lem:recsubres}.  Furthermore, for
  $k=2,\ldots,t-1$, define $\tilde{R}_k=(\tilde{R}_{k-1})^{b_k}B_k$
  with $\tilde{R}_0=\tilde{R}_1=1$.  Then, we have
  \begin{equation}
    \label{eq:nessubreslem}
    \nessubres_{k,j}(F,G)=
    (\tilde{R}_{k-1})^{b_{k,j}}\cdot\subres_j(P_1^{(k)},P_2^{(k)}).
  \end{equation}
\end{lem}
\begin{pf}
  By induction on $k$. For $k=1$, it is obvious by the definition of
  the nested subresultant.  Assume that (\ref{eq:nessubreslem}) is
  valid for $1,\ldots,k-1$.  Then, by the Fundamental Theorem of
  subresultants (Theorem~\ref{th:fundsubres}), we have
  \begin{equation*}
    \begin{split}
      \nessubres_{k-1,j_{k-1}}(F,G)
      &= (\tilde{R}_{k-2})^{b_{k-1}}B_{k-1}P_{l_{k-1}}^{(k-1)}
      = (\tilde{R}_{k-1})P_{1}^{(k)}, \\
      \frac{d}{dx}
      \left(
        \nessubres_{k-1,j_{k-1}}(F,G)
      \right)
      &= (\tilde{R}_{k-1})\frac{d}{dx}
      \left(
        P_1^{(k)}
      \right)
      = (\tilde{R}_{k-1})P_2^{(k)}.
    \end{split}
  \end{equation*}
  Then, we have
  \begin{equation}
    \label{eq:lem-nessubres-subres}
    \begin{split}
      \tilde{N}^{(k,j)}(F,G) &=
      (\tilde{R}_{k-1})\ N^{(j)}(P_1^{(k)},P_2^{(k)}),
      \\
      |\tilde{N}_\tau^{(k,j)}(F,G)| &=
      (\tilde{R}_{k-1})^{b_{k,j}} \ |N_\tau^{(j)}(P_1^{(k)},P_2^{(k)})|,
    \end{split}
  \end{equation}
  for $\tau=j,\ldots,0$.  Therefore, we have
  (\ref{eq:nessubreslem}), which proves the lemma.  \qed
\end{pf}
\begin{thm}
  \label{th:res-nes-subres-equiv}
  Let $F$ and $G$ be defined as in \textup{(\ref{eq:fg})}, and let
  $(P_1^{(1)},\ldots,P_{l_1}^{(1)},\ldots,$
  $P_1^{(t)},\ldots,P_{l_t}^{(t)})$ be complete recursive PRS for $F$
  and $G$ as in Definition~\ref{def:recprs}.  For $k=2,\ldots,t$ and
  $j=j_{k-1}-2,\ldots,0$, define $u_{k,j}$, $b_{k,j}$, $r_{k,j}$ as in
  Lemma~\ref{lem:recsubres} and $R'_k=(R'_{k-1})^{b_k}r_k$ with
  $R'_0=R'_1=1$.  Then, we have
  \begin{equation*}
    \recsubres_{k,j}(F,G)=(R'_{k-1})^{b_{k,j}}r_{k,j}\cdot\nessubres_{k,j}(F,G).
  \end{equation*}
\end{thm}
\begin{pf}
  By induction on $k$.  For $k=1$, it is obvious by the definitions of
  the recursive and the nested subresultants.  We first show that 
   $\bar{R}_k = \tilde{R}_k\cdot R'_k$
  for $k=0,\ldots,t-1$.  It is obvious for $k=0$ and $1$.  Let us assume
  $\bar{R}_{k-1} = \tilde{R}_{k-1}\cdot R'_{k-1}$. Then, we have
  \begin{equation*}
    \begin{split}
      \bar{R}_k &= (\bar{R}_{k-1})^{b_k} r_k B_k
      = (\tilde{R}_{k-1} \cdot R'_{k-1})^{b_k} r_k B_k
      = (\tilde{R}_{k-1})^{b_k}B_k\cdot (R'_{k-1})^{b_k}r_k \\
      &= \tilde{R}_k\cdot R'_k.
    \end{split}
  \end{equation*}
  Now, by Lemma~\ref{lem:recsubres}, we have
    $\recsubres_{k,j}(F,G) =
    (\bar{R}_{k-1})^{b_{k,j}}r_{k,j}\cdot
    \subres_j(P_1^{(k)},P_2^{(k)})$,
  then, by Lemma~\ref{lem:nessubres}, we have
  \begin{equation*}
    \begin{split}
      \recsubres_{k,j}(F,G)
      &= (\tilde{R}_{k-1}\cdot R'_{k-1})^{b_{k,j}}r_{k,j} \cdot
      \subres_j(P_1^{(k)},P_2^{(k)})  
      \\
      &= (R'_{k-1})^{b_{k,j}}r_{k,j} \cdot
      (\tilde{R}_{k-1})^{b_{k,j}} \cdot \subres_j(P_1^{(k)},P_2^{(k)})  
      \\
      &= (R'_{k-1})^{b_{k,j}}r_{k,j} \cdot \nessubres_{k,j}(F,G),
    \end{split}
  \end{equation*}
  which proves the theorem.
  \qed
\end{pf}
\begin{rem}
  Since $r_{k,j}=\pm1$, we see that $R'_k=\pm1$ hence the nested
  subresultant is equal to the recursive subresultant up to a sign.
  \qed
\end{rem}

\section{Reduced Nested Subresultants}
\label{sec:rednessubres}

The nested subresultant matrix has ``nested'' representation of
subresultant matrices, which makes practical use difficult.  However,
in some cases, we can reduce the representation of the nested
subresultant matrix to a ``flat'' representation, or a representation
without nested determinants by the Gaussian elimination; this is the
reduced nested subresultant (matrix).  As we will see, the size of the
reduced nested subresultant matrix becomes much smaller than that of
the recursive subresultant matrix, with reasonable computing time.

First, we illustrate the idea of reduction of the nested subresultant
matrix with an example.
\begin{exmp}
  \label{exp:rednessubres-1}
  Let $F(x)$ and $G(x)$ be defined as
  \begin{equation*}
    \begin{split}
      F(x) &= a_6 x^6+a_5 x^5+\cdots+a_0,\quad a_6\ne 0,\\
      G(x) &= b_5 x^5+b_4 x^4+\cdots+b_0,\quad b_5\ne 0,
    \end{split}
  \end{equation*}
  with vectors of coefficients $(a_6,a_5)$ and $(b_5,b_4)$ are
  linearly independent as vectors over $K$.  Assume that
  $\prs(F,G)=(P_1^{(1)}=F,\ P_2^{(1)}=G,\ P_3^{(1)}=\gcd(F,G))$ with
  $\deg(P_3^{(1)})=4$.  Consider the $(2,2)$-th nested subresultant;
  its matrix $\tilde{N}^{(2,2)}(F,G)$ is defined as
  \begin{equation}
    \label{eq:aj}
    \tilde{N}^{(2,2)}(F,G)=
    \begin{pmatrix}
      A_4 & 4A_4 & \\
      A_3 & 3A_3 & 4A_4 \\
      A_2 & 2A_2 & 3A_3 \\
      A_1 & A_1  & 2A_2 \\
      A_0 &      & A_1
    \end{pmatrix}
    ,
    \quad
    A_j=
    \begin{vmatrix}
      a_6 & b_5 & \\
      a_5 & b_4 & b_5 \\
      a_j & b_{j-1} & b_j
    \end{vmatrix}
    ,
  \end{equation}
  for $j\le 4$ with $b_j=0$ for $j<0$.  Now, let us calculate the
  leading coefficient of $\tilde{\subres}_{2,2}(F,G)$ as
  \begin{equation}
    \label{eq:ex-h1}
    \begin{split}
    |\tilde{N}^{(2,2)}_{2}|=
    \begin{vmatrix}
      A_4 & 4A_4 & \\
      A_3 & 3A_3 & 4A_4 \\
      A_2 & 2A_2 & 3A_3
    \end{vmatrix}
    &=
      \begin{vmatrix}
        \begin{vmatrix}
          a_6 & b_5 \\
          a_5 & b_4 & b_5 \\
          a_4 & b_3 & b_4
        \end{vmatrix}
        &
        \begin{vmatrix}
          a_6 & b_5 & \\
          a_5 & b_4 & b_5 \\
          4a_4 & 4b_3 & 4b_4
        \end{vmatrix}
        &
        \begin{vmatrix}
          a_6 & b_5 & \\
          a_5 & b_4 & b_5 \\
          0a_4 & 0b_3 & 0b_4
        \end{vmatrix}
        \\
        \noalign{\vskip2pt}
        \begin{vmatrix}
          a_6 & b_5 \\
          a_5 & b_4 & b_5 \\
          a_3 & b_2 & b_3
        \end{vmatrix}
        &
        \begin{vmatrix}
          a_6 & b_5 \\
          a_5 & b_4 & b_5 \\
          3a_3 & 3b_2 & 3b_3
        \end{vmatrix}
        &
        \begin{vmatrix}
          a_6 & b_5 & \\
          a_5 & b_4 & b_5 \\
          4a_4 & 4b_3 & 4b_4
        \end{vmatrix}
        \\
        \noalign{\vskip2pt}
        \noalign{\vskip2pt}
        \begin{vmatrix}
          a_6 & b_5 \\
          a_5 & b_4 & b_5 \\
          a_2 & b_1 & b_2
        \end{vmatrix}
        &
        \begin{vmatrix}
          a_6 & b_5 \\
          a_5 & b_4 & b_5 \\
          2a_2 & 2b_1 & 2b_2
        \end{vmatrix}
        &
        \begin{vmatrix}
          a_6 & b_5 & \\
          a_5 & b_4 & b_5 \\
          3a_3 & 3b_2 & 3b_3
        \end{vmatrix}
      \end{vmatrix}
      \\
    &=|H|
    =
    \left|
      \begin{pmatrix}
        H_{p,q}
      \end{pmatrix}
    \right|
    .
    \end{split}
  \end{equation}
  Then, we make the $(3,1)$ and the $(3,2)$ elements in $H_{p,q}$
  ($p,q=1,2,3$) equal to $0$ by adding the first and the second rows,
  multiplied by certain numbers, to the third row.  For example, in
  $H_{1,1}$, calculate $x_{11}$ and $y_{11}$ by solving a system of
  linear equations
  \begin{equation}
    \label{eq:h11}
    \left\{
      \begin{split}
        a_6 x_{11} + a_5 y_{11} &= -a_4 \\
        b_5 x_{11} + b_4 y_{11} &= -b_3
      \end{split}
    \right.
    ,
  \end{equation}
  (Note that (\ref{eq:h11}) has a solution in $K$ by the assumption),
  followed by adding the first row multiplied by $x_{11}$ and the
  second row multiplied by $y_{11}$, respectively, to the third row.
  Then, we have
  \begin{equation}
    \label{eq:H11}
    H_{1,1}=
    \begin{vmatrix}
      a_6 & b_5 \\
      a_5 & b_4 & b_5 \\
      0   & 0   & h_{11}
    \end{vmatrix}
    =
    \begin{vmatrix}
      a_6 & b_5 \\
      a_5 & b_4 
    \end{vmatrix}
    h_{11}
    \quad
    \mbox{with}
    \quad
    h_{11} = b_4 + y_{11}b_5.
  \end{equation}
  Doing similar calculations for the other $H_{p,q}$, we calculate
  $h_{p,q}$ ($p,q=1,2,3$) for $H_{p,q}$ similarly as in
  (\ref{eq:H11}).  Finally, by putting such new representations of
  $H_{p,q}$ into (\ref{eq:ex-h1}), we have
  \begin{equation}
    \label{eq:ex-nhat}
    |\tilde{N}^{(2,2)}_2|=
    \begin{vmatrix}
      a_6 & b_5 \\
      a_5 & b_4
    \end{vmatrix}
    ^3
    \begin{vmatrix}
      h_{11} & h_{12} & h_{13}\\
      h_{21} & h_{22} & h_{23}\\
      h_{31} & h_{32} & h_{33}
    \end{vmatrix}
    =
    \begin{vmatrix}
      a_6 & b_5 \\
      a_5 & b_4
    \end{vmatrix}
    ^3
    |\hat{N}^{(2,2)}_2|
    ,
  \end{equation}
  note that we have derived $\hat{N}^{(2,2)}_2$ as a reduced form of
  $\tilde{N}^{(2,2)}_2$.
  \qed
\end{exmp}

As (\ref{eq:ex-nhat}) shows, we derive a ``reduced'' form of the
nested subresultant matrix by the Gaussian elimination for solving
certain systems of linear equations.  We define the reduced nested
subresultant (matrix), as follows.

\begin{defn}[Reduced Nested Subresultant Matrix]
  \label{def:rednessubresmat}
  Let $F$ and $G$ \linebreak be defined as in \textup{(\ref{eq:fg})},
  and let $(P_1^{(1)},\ldots,P_{l_1}^{(1)},\ldots,$
  $P_1^{(t)},\ldots,P_{l_t}^{(t)})$ be complete recursive PRS for $F$
  and $G$ as in Definition~\ref{def:recprs}.  Then, for each pair of
  numbers $(k,j)$ with $k=1,\ldots,t$ and $j=j_{k-1}-2,\ldots,0$,
  define matrix $\hat{N}^{(k,j)}(F,G)$ recursively as follows.
  \begin{enumerate}
  \item For $k=1$, let $\hat{N}^{(1,j)}(F,G)=N^{(j)}(F,G)$.
  \item For $k>1$, let $\hat{N}_U^{(k-1,j_{k-1})}(F,G)$ be a
    sub-matrix of $\hat{N}^{(k-1,j_{k-1})}(F,G)$ by deleting the
    bottom $j_{k-1}+1$ rows, and $\hat{N}_L^{(k-1,j_{k-1})}(F,G)$ be a
    sub-matrix of \linebreak $\hat{N}^{(k-1,j_{k-1})}(F,G)$ by taking
    the bottom $j_{k-1}+1$ rows, respectively.  For
    $\tau=j_{k-1},\ldots,0$ let $\hat{N}_\tau^{(k-1,j_{k-1})}(F,G)$ be
    a sub-matrix of $\hat{N}^{(k-1,j_{k-1})}(F,G)$ by putting
    $\hat{N}_U^{(k-1,j_{k-1})}(F,G)$ on the top and the
    $(j_{k-1}-\tau+1)$-th row of $\hat{N}_L^{(k-1,j_{k-1})}(F,G)$ in
    the bottom row.  Let
    $\hat{A}_\tau^{(k-1)}=|\hat{N}_\tau^{(k-1,j_{k-1})}|$ and
    construct a matrix $H^{(k,j)}$ as
  \begin{equation}
    \label{eq:h1}
    H^{(k,j)} =
    \begin{pmatrix}
      H^{(k,j)}_{p,q}
    \end{pmatrix}
    =
    N^{(j)}
    \left(
      \hat{A}^{(k-1)}(x), \frac{d}{dx}\hat{A}^{(k-1)}(x)
    \right),
  \end{equation}
  where
  \begin{equation*}
    \hat{A}^{(k-1)}(x)=\hat{A}^{(k-1)}_{j_{k-1}}x^{j_{k-1}}+\cdots
    +\hat{A}^{(k-1)}_0x^0.
  \end{equation*}
  Since $\hat{N}_\tau^{(k-1,j_{k-1})}$ consists of
  $\hat{N}_U^{(k-1,j_{k-1})}$ and a row vector in the bottom, we
  express $\hat{N}_U^{(k-1,j_{k-1})}=
  \begin{pmatrix}
    U^{(k)} | \bm{v}^{(k)}
  \end{pmatrix}
  $, where $U^{(k)}$ is a square matrix and $\bm{v}^{(k)}$ is a column
  vector, and the row vector in the bottom by $
  \begin{pmatrix}
    \bm{b}_{p,q}^{(k,j)} \Bigm| g_{p,q}^{(k,j)}
  \end{pmatrix}
  $, where $\bm{b}_{p,q}^{(k,j)}$ is a row vector and
  $g_{p,q}^{(k,j)}$ is a number, respectively, such that
  \begin{equation}
    \label{eq:Hpq}
    H_{p,q}^{(k,j)}=
    \left|
      \begin{array}{c|c}
        U^{(k)} & \bm{v}^{(k)} \\
        \hline
        \bm{b}_{p,q}^{(k,j)} & g_{p,q}^{(k,j)}
      \end{array}
    \right|
    ,
  \end{equation}
  with $\bm{b}_{p,q}^{(k,j)}=\bm{0}$ and $g_{p,q}^{(k,j)}=0$ for
  $H_{p,q}^{(k,j)}=0$.  Furthermore, we assume that $U^{(k)}$ is not
  singular.  Then, for $p=1,\ldots,n_1^{(k)}+n_2^{(k)}-j$ and
  $q=2,\ldots,n_1^{(k)}+n_2^{(k)}-j$, calculate a row vector
  $\bm{x}_{p,q}^{(k)}$ as a solution of the equation
  \begin{equation}
    \label{eq:ukeq}
    \bm{x}_{p,q}^{(k,j)}U^{(k)}=- \bm{b}_{p,q}^{(k,j)},
  \end{equation}
  and define $h_{p,q}^{(k,j)}$ as
  \begin{equation}
    \label{eq:hpq}
    h_{p,q}^{(k,j)}=\bm{x}_{p,q}^{(k,j)}\bm{v}^{(k,j)}.
  \end{equation}
  Note that we have
  \begin{equation}
    \label{eq:Hpq-2}
    H_{p,q}^{(k,j)}=
    \left|
      \begin{array}{c|c}
        U^{(k)} & \bm{v}^{(k)} \\
        \hline
        \bm{0} & h_{p,q}^{(k,j)}
      \end{array}
    \right|
    =
    \left|
      U^{(k)}
    \right|
    h_{p,q}^{(k,j)}
    .
  \end{equation}
  Finally, define $\hat{N}^{(k,j)}(F,G)$ as
  \begin{equation}
    \label{eq:rednessubresmat}
    \hat{N}^{(k,j)}(F,G)=
    \begin{pmatrix}
      h_{1,1}^{(k,j)} & h_{1,2}^{(k,j)} & \cdots & h_{1,J_{k,j}}^{(k,j)} \\
      h_{2,1}^{(k,j)} & h_{2,2}^{(k,j)} & \cdots & h_{2,J_{k,j}}^{(k,j)} \\
      \vdots & \vdots & & \vdots \\
      h_{I_{k,j},1}^{(k,j)} & h_{I_{k,j},2}^{(k,j)} & \cdots &
      h_{I_{k,j},J_{k,j}}^{(k,j)} 
    \end{pmatrix}
    ,
  \end{equation}
  where
  \begin{equation}
    \label{eq:ikj-jkj}
    \begin{split}
      I_{k,j} &= n_1^{(k)}+n_2^{(k)}-j=(2j_{k-1}-2j-1)+j,
      \\
      J_{k,j} &= n_1^{(k)}+n_2^{(k)}-2j=2j_{k-1}-2j-1.
    \end{split}
  \end{equation}
  Then, $\hat{N}^{(k,j)}(F,G)$ is called the \emph{$(k,j)$-th reduced
    nested subresultant matrix of $F$ and $G$}. \qed
  \end{enumerate}
\end{defn}
\begin{rem}
  \label{rem:rednessubresmatord}
  Definition~\ref{def:rednessubresmat} shows that, For $k=1,\ldots,t$
  and $j<j_{k-1}-1$, the numbers of rows and columns of the $(k,j)$-th
  reduced nested subresultant matrix $\hat{N}^{(k,j)}(F,G)$ are
  $I_{k,j}$ and $J_{k,j}$ in (\ref{eq:ikj-jkj}), respectively, which
  are much smaller than those of the recursive subresultant matrix of
  the corresponding degree (see Proposition~\ref{prop:recsubresmat}).
  \qed
\end{rem}
\begin{defn}[Reduced Nested Subresultant]
  \label{def:rednessubres}
  Let $F$ and $G$ be defined as in \textup{(\ref{eq:fg})}, and let
  $(P_1^{(1)},\ldots,P_{l_1}^{(1)},\ldots,P_1^{(t)},\ldots,P_{l_t}^{(t)})$
  be complete recursive PRS for $F$ and $G$ as in
  Definition~\ref{def:recprs}.  For $j=j_{k-1}-2,\ldots,0$ and
  $\tau=j,\ldots,0$, let
  $\hat{N}_\tau^{(k,j)}=\hat{N}_\tau^{(k,j)}(F,G)$ be a sub-matrix of
  the $(k,j)$-th reduced nested subresultant matrix
  $\hat{N}^{(k,j)}(F,G)$ obtained by the top
  $n_1^{(k)}+n_2^{(k)}-2j-1$ rows and the
  $(n_1^{(k)}+n_2^{(k)}-j-\tau)$-th row (note that
  $\hat{N}_\tau^{(k,j)}(F,G)$ is a square matrix).  Then, the
  polynomial
  \begin{equation*}
    \rednessubres_{k,j}(F,G)=|\hat{N}_j^{(k,j)}(F,G)|x^j+\cdots
    +|\hat{N}_0^{(k,j)}(F,G)|x^0
  \end{equation*}
  is called the $(k,j)$-th reduced nested subresultant of $F$ and $G$.
  \qed
\end{defn}

Now, we derive the relationship between the nested and the reduced
nested subresultants.

\begin{thm}
  \label{th:nessubres-rednessubres-equiv}
  Let $F$ and $G$ be defined as in \textup{(\ref{eq:fg})}, and let
  $(P_1^{(1)},\ldots,P_{l_1}^{(1)},\ldots,$
  $P_1^{(t)},\ldots,P_{l_t}^{(t)})$ be complete recursive PRS for $F$
  and $G$ as in Definition~\ref{def:recprs}.  For $k=2,\ldots,t$,
  $j=j_{k-1}-2,\ldots,0$ with $J_{k,j}$ as in
  (\ref{eq:rednessubresmat}), define $\hat{B}_{k,j}$ and $\hat{R}_k$
  as
  \begin{equation*}
    \hat{B}_{k,j}=|U^{(k)}|^{J_{k,j}}
  \end{equation*}
  with $\hat{B}_k=\hat{B}_{k,j_k}$ and $\hat{B}_1=\hat{B}_2=1$, and
  \begin{equation*}
    \hat{R}_k = (\hat{R}_{k-1}\cdot\hat{B}_{k-1})^{J_{k,j_k}}
  \end{equation*}
  with $\hat{R}_1=\hat{R}_2=1$, respectively.  Then, we have
  \begin{equation*}
    \nessubres_{k,j}(F,G)=
    (\hat{R}_{k-1}\cdot\hat{B}_{k-1})^{J_{k,j}}\hat{B}_{k,j}\cdot
    \rednessubres_{k,j}(F,G).
  \end{equation*}
\end{thm}
To prove Theorem~\ref{th:nessubres-rednessubres-equiv}, we prove the
following lemma.
\begin{lem}
  \label{le:nessubres-rednessubres-equiv}
  For $k=1,\ldots,t$, $j=j_{k-1}-2,\ldots,0$ and $\tau=j,\ldots,0$, we
  have
  \begin{equation}
    \label{eq:nessubresmat-rednessubresmat-equiv}
    |\tilde{N}_\tau^{(k,j)}(F,G)|=
    (\hat{R}_{k-1}\cdot\hat{B}_{k-1})^{J_{k,j}}\hat{B}_{k,j}
    |\hat{N}_\tau^{(k,j)}(F,G)|.
  \end{equation}
\end{lem}
\begin{pf}
  By induction on $k$.  For $k=1$, it is obvious from the definitions
  of the nested and the reduced nested subresultants.  Assume that
  (\ref{eq:nessubresmat-rednessubresmat-equiv}) is valid for
  $1,\ldots,k-1$.  Then, for $\tau=j_{k-1},\ldots,0$, we have
  \begin{equation*}
    \begin{split}
      |\tilde{N}_\tau^{(k-1,j_{k-1})}(F,G)|
    =&
    (\hat{R}_{k-2}\cdot\hat{B}_{k-2})^{J_{k-1,j_{k-1}}}\hat{B}_{k-1,j_{k-1}}
    |\hat{N}_\tau^{(k-1,j_{k-1})}(F,G)|
    \\
    =&
    (\hat{R}_{k-1}\cdot\hat{B}_{k-1})|\hat{N}_\tau^{(k-1,j_{k-1})}(F,G)|.
    \end{split}
  \end{equation*}
  Let 
  \begin{equation*}
    \tilde{A}_\tau^{(k-1)}=|\tilde{N}_\tau^{(k-1,j_{k-1})}(F,G)|,
    \quad
    \hat{A}_\tau^{(k-1)}=|\hat{N}_\tau^{(k-1,j_{k-1})}(F,G)|,
  \end{equation*}
  and $H_\tau^{(k,j)}=
  \begin{pmatrix}
    H^{(k,j)}_{\tau_{p,q}}
  \end{pmatrix}
  $ be a sub-matrix of $H^{(k,j)}$ in (\ref{eq:h1}) by taking the top
  $J_{k,j}$ rows and the $(I_{k,j}-\tau)$-th row, where $I_{k,j}$ and
  $J_{k,j}$ are defined as in (\ref{eq:ikj-jkj}), respectively.  Then,
  we have
  \begin{equation*}
    \begin{split}
    & |H^{(k,j)}_\tau|
    \\
    =&
    {\small
    \begin{vmatrix}
      \hat{A}_{j_{k-1}}^{(k-1)} &        &       &
      j_{k-1}\hat{A}_{j_{k-1}}^{(k-1)}    &        & \\
      \vdots & \ddots &  & \vdots & \ddots & \\
      \vdots &        & \hat{A}_{j_{k-1}}^{(k-1)} &
      \vdots &        & j_{k-1}\hat{A}_{j_{k-1}}^{(k-1)}\\
      \vdots &        & \vdots & \vdots &        & \vdots \\
      \hat{A}_{2j-j_{k-1}+3}^{(k-1)} & \cdots &
      \hat{A}_{j+1}^{(k-1)} &
      (2j-j_{k-1}+3)\hat{A}_{2j-j_{k-1}+3}^{(k-1)} & \cdots &
      (j+2)\hat{A}_{j+2}^{(k-1)} \\
      \hline
      \hat{A}_{j-j_{k-1}+\tau+2}^{(k-1)} & \cdots &
      \hat{A}_\tau^{(k-1)} &
      (j-j_{k-1}+\tau+2)\hat{A}_{j-j_{k-1}+\tau+2}^{(k-1)} &
      \cdots & (\tau+1)\hat{A}_{\tau+1}^{(k-1)}  
    \end{vmatrix}
    ,
  }
  \end{split}
  \end{equation*}
  where $\hat{A}_{l}^{(k-1)}=0$ for $l<0$.  On the other hand, by the
  definition of the $(k,j)$-th nested subresultant, we have
  \begin{align}
    & |\tilde{N}_\tau^{(k,j)}(F,G)| \nonumber
    \\
    =&
    \small{
    \begin{vmatrix}
      \tilde{A}_{j_{k-1}}^{(k-1)} &        &       &
      j_{k-1}\tilde{A}_{j_{k-1}}^{(k-1)}    &        & \\
      \vdots & \ddots &  & \vdots & \ddots & \\
      \vdots &        & \tilde{A}_{j_{k-1}}^{(k-1)} &
      \vdots &        & j_{k-1}\tilde{A}_{j_{k-1}}^{(k-1)}\\
      \vdots &        & \vdots & \vdots &        & \vdots \\
      \tilde{A}_{2j-j_{k-1}+3}^{(k-1)} & \cdots &
      \tilde{A}_{j+1}^{(k-1)} &
      (2j-j_{k-1}+3)\tilde{A}_{2j-j_{k-1}+3}^{(k-1)} & \cdots &
      (j+2)\tilde{A}_{j+2}^{(k-1)} \\
      \hline
      \tilde{A}_{j-j_{k-1}+\tau+2}^{(k-1)} & \cdots &
      \tilde{A}_\tau^{(k-1)} &
      (j-j_{k-1}+\tau+2)\tilde{A}_{j-j_{k-1}+\tau+2}^{(k-1)} &
      \cdots & (\tau+1)\tilde{A}_{\tau+1}^{(k-1)}  
    \end{vmatrix}
  }
  \nonumber
    \\
    \label{eq:nesresmat-tau-2}
    =&
    (\hat{R}_{k-1}\cdot\hat{B}_{k-1})^{J_{k,j}}|H^{(k,j)}_\tau|,
  \end{align}
  where $\tilde{A}_{l}^{(k-1)}=0$ for $l<0$.
  (Note that $\tilde{N}_\tau^{(k,j)}$ and $H^{(k,j)}_\tau$ are square
  matrices of order $J_{k,j}$.)  By 
  Definition~\ref{def:rednessubresmat}, we can express
  $H^{(k,j)}_{\tau_{p,q}}$ as
  \begin{equation*}
    H^{(k,j)}_{\tau_{p,q}}=
    \left|
      \begin{array}{c|c}
        U^{(k)} & \bm{v}^{(k)} \\
        \hline
        \bm{b}'^{(k,j)}_{p,q} & g'^{(k,j)}_{p,q}
      \end{array}
    \right|
  ,
  \end{equation*}
  with $\bm{b}'^{(k,j)}_{p,q}=\bm{0}$ and $g'^{(k,j)}_{p,q}=0$ for
  $H^{(k,j)}_{\tau_{p,q}}=0$.  Note that, for $q=1,\ldots,J_{k,j}$, we have
  $\bm{b}'^{(k,j)}_{p,q}=\bm{b}^{(k,j)}_{p,q}$ and
  $g'^{(k,j)}_{p,q}=g^{(k,j)}_{p,q}$ for $p=1,\ldots,J_{k,j}-1$, and
  $\bm{b}'^{(k,j)}_{J_{k,j},q}=\bm{b}^{(k,j)}_{I_{k,j}-\tau,q}$ and
  $g'^{(k,j)}_{J_{k,j},q}=g^{(k,j)}_{I_{k,j}-\tau,q}$, where
  $\bm{b}^{(k,j)}_{p,q}$ and $g^{(k,j)}_{p,q}$ are defined as in
  (\ref{eq:Hpq}), respectively.
  Thus, by (\ref{eq:ukeq})--(\ref{eq:rednessubresmat}), we have
  \begin{equation}
    \label{eq:nesresmat-tau-2-reduced}
    |H^{(k,j)}_\tau|=|U^{(k)}|^{J_{k,j}}|\hat{N}^{(k,j)}_\tau(F,G)|
    =\hat{B}_{k,j}|\hat{N}^{(k,j)}_\tau(F,G)|,
  \end{equation}
  and, by putting (\ref{eq:nesresmat-tau-2-reduced}) into
  (\ref{eq:nesresmat-tau-2}), we prove the lemma.
  \qed
\end{pf}
\begin{rem}
  \label{re:nessubresmat}
  We can calculate the $(k,j)$-th reduced nested subresultant matrix
  as a sub-matrix of the $(k,0)$-th reduced nested subresultant
  matrix.  In~\textup{(\ref{eq:h1})}, we see that the matrix
  $H^{(k,j)}$ is a sub-matrix of %
  $N
  \left( \hat{A}^{(k-1)}(x),
    \frac{d}{dx}\hat{A}^{(k-1)}(x)
  \right)
  $, %
  and, by the construction of the reduced nested resultant matrix
  \textup{(\ref{eq:rednessubresmat})}, we see that
  $\hat{N}^{(k,j)}(F,G)$ is a sub-matrix of $\hat{N}^{(k,0)}(F,G)$ by
  taking the left $n_2^{(k)}-j$ columns from those corresponding to
  the coefficients of $\hat{A}^{(k-1)}(x)$ and the left
  $n_1^{(k)}-j$ columns from those corresponding to the coefficients of
  $\frac{d}{dx}\hat{A}^{(k-1)}(x)$, then taking the top
  $n_1^{(k)}+n_2^{(k)}-j$ rows. 
  \qed
\end{rem}
\begin{rem}
  \label{re:time}
  We can estimate arithmetic computing time for the $(k,j)$-th reduced
  nested resultant matrix $\hat{N}^{(k,j)}$ in
  \textup{(\ref{eq:rednessubresmat})}, as follows.  The computing time
  for the elements $h_{p,q}$ is dominated by the time for the Gaussian
  elimination of $U^{(k)}$.  Since the order of $U^{(k)}$
  ($k=2,\ldots,t$) is equal to $2(j_{k-2}-j_{k-1}-1)$ (see
  Remark~\ref{rem:rednessubresmatord}), it is bounded by
  $O((j_{k-2}-j_{k-1})^3)$ (see Golub and van Loan
  \textup{\cite{gol-vloa1996}}).  As Remark~\ref{re:nessubresmat}
  shows, we can calculate $\hat{N}^{(k,j)}(F,G)$ for $j<j_{k-1}-2$ by
  $\hat{N}^{(k,0)}(F,G)$.  Therefore, total computing time for
  $\hat{N}^{(k,j)}$ for entire recursive PRS ($k=1,\ldots,t$) is
  bounded by
  \begin{equation*}
    \sum_{k=2}^t O((j_{k-2}-j_{k-1})^3)=
    O
    \left(
      \sum_{k=2}^t (j_{k-2}-j_{k-1})^3
    \right)
    =
    O((j_0-j_{t-1})^3)
    =O(m^3),
  \end{equation*}
  note that $j_0=m$ (see Remark~\ref{re:recprs}).  See also for the
  concluding remarks. \qed
\end{rem}

\section{Concluding Remarks}
\label{sec:disc}

In this paper, we have introduced concepts of \emph{recursive PRS} and
\emph{recursive subresultants}, and investigated constructions of
their subresultant matrices to compute the recursive subresultants.
Among three different constructions of recursive subresultant
matrices, we have shown that the reduced nested subresultant matrix
reduces the size of the matrix drastically to at most the order of the
degree of initial polynomials in each PRSs, compared with the naive
recursive subresultant matrix.  We have also shown that we can
calculate the reduced nested subresultant matrix by the Gaussian
elimination of order at most the sum of the degree of
initial polynomials in each PRSs.

From a point of view of computational complexity, the algorithm for
the reduced nested subresultant matrix has a cubic complexity bound in
terms of the degree of the input polynomials (see
Remark~\ref{re:time}).  However, subresultant algorithms which have a
quadratic complexity bound in terms of the degree of the input
polynomials have been proposed (\cite{duc2000},
\cite{lom-roy-eldin2000}); those algorithms exploit the structure of the
Sylvester matrix to increase their efficiency with controlling the
size of coefficients well.  Although, in this paper, we have primarily
focused our attention into reducing the structure of the nested
subresultant matrix to ``flat'' representation, development of more
efficient algorithms such as exploiting the structure of the Sylvester
matrix would be the next problem.  Furthermore, the reduced nested
subresultant may involve fractions which may be unusual for
subresultants, thus more detailed analysis of computational
efficiency including comparison with (ordinary and recursive)
subresultants would also be necessary.

We expect that the reduced nested subresultants can be used for
approximate algebraic computation such as the square-free
decomposition of approximate univariate polynomials with approximate
GCD computations based on Singular Value Decomposition (SVD) of
subresultant matrices (\cite{cor-gia-tra-wat1995},
\cite{emi-gal-lom1997}), which motivates the present work.  For the
approximate square-free decomposition of the given polynomial $P(x)$,
we have to calculate the approximate GCDs of $P(x),\ldots,$
$P^{(n)}(x)$ (by $P^{(n)}(x)$ we denote the $n$-th derivative of
$P(x)$) or those of the recursive PRS for $P(x)$ and $P'(x)$; we have
to find the representation of the subresultant matrices for
$P(x),\ldots,P^{(n)}(x)$, or that for the recursive PRS for $P(x)$ and
$P'(x)$, respectively.  As for the former approach, several algorithms
based on different representations of subresultant matrices have been
proposed (\cite{dtoc-gveg2002}, \cite{rup1999}); our reduced nested
subresultant matrix can be used as for the latter approach.  To make
use of the reduced nested subresultant matrix, we need to reveal the
relationship between the structure of the subresultant matrices and
their singular values; this is the problem on which we are working
now.

\ack{The author would like to thank Prof.\ Tateaki Sasaki and
  anonymous referees for their valuable comments and suggestions that
  helped improve the paper.}

\newpage

\begin{figure}[t]
  \centering
  \begin{multline*}
    \bar{N}^{(k,j)}(F,G)=
    \\
    \tiny
  \left(
  \begin{array}{c|c|c|c|c|c|c|c}
    \hline
    \bar{N}_U^{(k-1,j_{k-1})} &            &        &          &
                  &            &        &         
     \\
     \cline{1-2}
             & \bar{N}_U^{(k-1,j_{k-1})}   &        &          &
             &            &        &          
     \\
     \cline{2-2}
             &            & \ddots &          &
             &            &        &          
     \\
     \cline{4-4}
             &            &        & \bar{N}_U^{(k-1,j_{k-1})} &
             &            &        &          
     \\
     \cline{4-5}
              &            &        &          &
     \bar{N}_U^{(k-1,j_{k-1})} &            &        &          
     \\
     \cline{5-6}
              &            &        &          &
              & \bar{N}_U^{(k-1,j_{k-1})}   &        &
     \\
     \cline{6-6}
              &            &        &          &
               &            & \ddots &         
     \\
     \cline{8-8}
               &            &        &         &
              &            &        & \bar{N}_U^{(k-1,j_{k-1})}
     \\
     \hline
              & 0 \cdots\cdots 0 &        &          &
              & 0 \cdots\cdots 0 &        &          
     \\
     \cline{2-2}\cline{6-6}
     \bar{N}_L^{(k-1,j_{k-1})} &            &        &          &
     \bar{N}_L'^{(k-1,j_{k-1})} &            &        & 
     \\
              & \bar{N}_L^{(k-1,j_{k-1})}   &        &     &
              & \bar{N}_L'^{(k-1,j_{k-1})}   &        &          
    \\       
    \cline{1-1}\cline{5-5}
               &            & \cdots &          &
              &            & \cdots &          
    \\
    \cline{2-2}\cline{4-4}\cline{6-6}\cline{8-8}
               &            &        &          &
              &            &        &          
    \\
              &            &        & \bar{N}_L^{(k-1,j_{k-1})} &
              &            &        & \bar{N}_L'^{(k-1,j_{k-1})}
    \\
              &            &        &          &
              &            &        &          
    \\
    \hline
  \end{array}
  \right).
  \end{multline*}
  \caption{Illustration of $\bar{N}^{(k,j)}(F,G)$.  Note that
    the number of blocks $\bar{N}_L^{(k-1,j_{k-1})}$ is $j_{k-1}-j-1$,
    whereas that of $\bar{N}_L^{'(k-1,j_{k-1})}$ is $j_{k-1}-j$; see
    Definition~\ref{def:recsubresmat} for details.}
  \label{fig:recsubresmat}
\end{figure}

\begin{figure}[b]
  \centering
  \begin{multline*}
    M^{(k,j)}(F,G)=
    \\
    \small
    \left(
      \begin{array}{c|c|c|c|c|c|c|c|c|c}
        \hline
        W_{k-1} & \bm{0}      &        &          &
        &            &        &         &
        \\
        \cline{1-2}
        *       & N_U^{(j_{k-1})}  &        &          &
        &            &        &         &
        \\
        \cline{1-2}
        &            & \ddots &          &
        &            &        &          &
        \\
        \cline{4-5}
        &            &        & W_{k-1}  &       
        \small
        \bm{0} &           &            &        &
        \\
        \cline{4-5}
        &            &        &  *      &
        N_U^{(j_{k-1})} &  &  &  &
        \\
        \cline{4-7}
        &            &        &          &       &
        \small
        W_{k-1}  &  \bm{0}    &     &
        \\
        \cline{6-7}
        &            &        &          &       &
        \small
        * & N_U^{(j_{k-1})}   &     &
        \\
        \cline{6-7}
        &            &        &          &       &
        &            & \ddots &        
        \\
        \cline{9-10}
        &            &        &         &
        &            &        &         &
        \small
        W_{k-1} & \bm{0}\\
        \cline{9-10}
        &            &        &         &
        &            &        &         &
        * & N_U^{(j_{k-1})}
        \\
        \hline
        * &  N_L^{(j_{k-1})}  &        &          &     &
        * & N_L'^{(j_{k-1})}
        &            &       
        \\
        \cline{1-2}\cline{6-7}
        &            & \ddots &          &    &
        &            & \ddots &          
        \\
        \cline{4-5}\cline{9-10}
        &            &        & * & N_L^{(j_{k-1})}
        &
        &            &        & 
        * & N_L'^{(j_{k-1})}
        \\
        \hline
      \end{array}
    \right).
  \end{multline*}
  \caption{Illustration of
    $M^{(k,j)}(F,G)$.  Note that the number of column blocks is
    equal to $b_{k,j}=2j_{k-1}-2j-1$; see Lemma~\ref{lem:recsubresmatm}
    for details.}
  \label{fig:recsubresmatm}
\end{figure}

\begin{figure}[htbp]
  \centering
  \begin{multline*}
    \bar{M}^{(k,j)}(F,G)=
    \\
    \left(
      \tiny
      \begin{array}{c|c|c|c|c|c|c|c|c|c}
        \hline
        \begin{array}{c|c}
          W_{k-1} & \bm{0}\\
          \hline
          * & \bar{N}_U^{(j_{k-1})}
        \end{array}
        & \bm{0} & & & & & & &
        \\
        \cline{1-2}
        & & \ddots & & & & & &
        \\
        \cline{4-5}
        & & &
        \begin{array}{c|c}
          W_{k-1} & \bm{0}\\
          \hline
          * & \bar{N}_U^{(j_{k-1})}
        \end{array}
        & \bm{0} & & & &
        \\
        \cline{4-7}
        & & & & & 
        \begin{array}{c|c}
          W_{k-1} & \bm{0}\\
          \hline
          * & \bar{N}_U^{(j_{k-1})}
        \end{array}
        & \bm{0} & & &
        \\
        \cline{6-7}
        & & & & & & & \ddots &
        \\
        \cline{9-10}
        & & & & & & & & 
        \begin{array}{c|c}
          W_{k-1} & \bm{0}\\
          \hline
          * & \bar{N}_U^{(j_{k-1})}
        \end{array}
        & \bm{0}
        \\
        \hline
        * & \bm{p}_1^{(k)} & & & & * & \bm{p}_2^{(k)} & & &
        \\
        \cline{1-2}\cline{6-7}
        & & \ddots & & & & & \ddots & &
        \\
        \cline{4-5}\cline{9-10}
        & & & * & \bm{p}_1^{(k)} & &  & & * & \bm{p}_2^{(k)}
        \\
        \hline
      \end{array}
    \right).
  \end{multline*}
  \caption{Illustration of
    $\bar{M}^{(k,j)}(F,G)$;  see Lemma~\ref{lem:recsubresmatmbar} for
    details.}
  \label{fig:recsubresmatmbar}
\end{figure}

\begin{figure}[htbp]
  \centering
  \begin{multline*}
    \hat{M}^{(k,j)}(F,G)=
    \\
    \tiny
    \left(
      \begin{array}{c|c|c|c|c|c|c|c|c|c|c|c}
        \hline
        \begin{array}{c|c}
          W_{k-1} & \bm{0}\\
          \hline
          * & \bar{N}_U^{(j_{k-1})}
        \end{array}
        & & & & &
        \\
        \cline{1-1}
        & \ddots & & & &
        \\
        \cline{3-3}
        & & 
        \begin{array}{c|c}
          W_{k-1} & \bm{0}\\
          \hline
          * & \bar{N}_U^{(j_{k-1})}
        \end{array}
        & & &
        \\
        \cline{3-4}
        & & &
        \begin{array}{c|c}
          W_{k-1} & \bm{0}\\
          \hline
          * & \bar{N}_U^{(j_{k-1})}
        \end{array}
        & &
        \\
        \cline{4-4}
        & & & & \ddots &
        \\
        \cline{6-6}
        & & & & &
        \begin{array}{c|c}
          W_{k-1} & \bm{0}\\
          \hline
          * & \bar{N}_U^{(j_{k-1})}
        \end{array}
        \\
        \hline
        * & & & * & & & \bm{p}_1^{(k)} & & & \bm{p}_2^{(k)} &
        \\
        \cline{1-1}\cline{4-4}
        & \ddots & & & \ddots & & & \ddots & & & \ddots &
        \\
        \cline{3-3}\cline{6-6}
        & & * & & & * & & & \bm{p}_1^{(k)} & & & \bm{p}_2^{(k)}
        \\
        \hline
      \end{array}
    \right)
  \end{multline*}
  \caption{Illustration of $\hat{M}^{(k,j)}(F,G)$.  Note that the
    lower-right block which consists of $\bm{p}_1^{(k)}$ and
    $\bm{p}_2^{(k)}$ is equal to $N^{(j)}(P_1^{(k)},P_2^{(k)})$ and
    $|W_{k-1}|=|\bar{N}_{U}^{(j_{k-1})}|=1$; see
    Lemma~\ref{lem:recsubresmat} for details.}
  \label{fig:recsubresmathat}
\end{figure}



\begin{thebibliography}{10}

\expandafter\ifx\csname url\endcsname\relax
  \def\url#1{\texttt{#1}}\fi
\expandafter\ifx\csname urlprefix\endcsname\relax\def\urlprefix{URL }\fi

\bibitem{boc-cos-roy1998}
J.~Bochnak, M.~Coste, M.-F. Roy, Real algebraic geometry, vol.~36 of Ergebnisse
  der Mathematik und ihrer Grenzgebiete (3) [Results in Mathematics and Related
  Areas (3)], Springer-Verlag, Berlin, 1998, translated from the 1987 French
  original, Revised by the authors.

\bibitem{bro-tra71}
W.~S. Brown, J.~F. Traub, {O}n {E}uclid's {A}lgorithm and the {T}heory of
  {S}ubresultants, J.~ACM 18~(4) (1971) 505--514.

\bibitem{col1967}
G.~E. Collins, {S}ubresultants and {R}educed {P}olynomial {R}emainder
  {S}equences, J. ACM 14~(1) (1967) 128--142.

\bibitem{cor-gia-tra-wat1995}
R.~M. Corless, P.~M. Gianni, B.~M. Trager, S.~M. Watt, The singular value
  decomposition for polynomial systems, in: Proc. ISSAC 1995, ACM, 1995.

\bibitem{dtoc-gveg2002}
G.~M. Diaz-Toca, L.~Gonzalez-Vega, Barnett's theorems about the greatest common
  divisor of several univariate polynomials through {B}ezout-like matrices, J.
  Symbolic Comput. 34~(1) (2002) 59--81.

\bibitem{duc2000}
L.~Ducos, Optimizations of the subresultant algorithm, J. Pure Appl. Algebra
  145~(2) (2000) 149--163.

\bibitem{emi-gal-lom1997}
I.~Z. Emiris, A.~Galligo, H.~Lombardi, Certified approximate univariate {GCD}s,
  J. Pure Appl. Algebra 117/118 (1997) 229--251, {A}lgorithms for algebra
  (Eindhoven, 1996).

\bibitem{vzg-luc2003}
J.~von~zur Gathen, T.~L{\"u}cking, Subresultants revisited, Theoret. Comput.
  Sci. 297~(1-3) (2003) 199--239, {L}atin {A}merican theoretical informatics
  (Punta del Este, 2000).

\bibitem{gol-vloa1996}
G.~H. Golub, C.~F. Van~Loan, Matrix Computations, 3rd ed., The Johns Hopkins
  University Press, 1996.

\bibitem{knu1998}
D.~Knuth, The Art of Computer Programming, vol. 2: Seminumerical Algorithms,
  3rd ed., Addison-Wesley, 1998.

\bibitem{lom-roy-eldin2000}
H.~Lombardi, M.-F. Roy, M.~S. El~Din, New structure theorem for subresultants,
  J. Symbolic Comput. 29~(4--5) (2000) 663--689.

\bibitem{loos1983}
R.~Loos, Generalized polynomial remainder sequences, in: B.~Buchberger, G.~E.
  Collins, R.~Loos (eds.), Computer Algebra: Symbolic and Algebraic
  Computation, 2nd ed., Springer-Verlag, 1983, pp. 115--137.

\bibitem{rup1999}
D.~Rupprecht, An algorithm for computing certified approximate {GCD} of $n$
  univariate polynomials, J. Pure Appl. Algebra 139 (1999) 255--284.

\bibitem{ter2003}
A.~Terui, Subresultants in recursive polynomial remainder sequence, in:
  V.~Ganzha, E.~Mayr, E.~Vorozhtsov (eds.), Proc.\ The 6th
  International Workshop on Computer Algebra in Scientific Computing:
  CASC 2003 (Passau, Germany, 2003), Institute f\"ur Informatik,
  Technische Universit\"at M\"unchen, Garching, 2003, pp. 363--375.

\bibitem{ter2005}
A.~Terui, Recursive polynomial remainder sequence and the nested subresultants,
  in: V.~Ganzha, E.~Mayr, E.~Vorozhtsov (eds.), Computer Algebra in Scientific
  Computing, vol. 3718 of Lecture Notes in Computer Science,
  Springer, 2005, pp. 445--456, CASC 2005 (Kalamata, Greece, 2005).
\end{thebibliography}
\end{document}